\documentclass[hidelinks,onefignum,onetabnum]{siamart250106}

\newcommand{\mycomment}[1]{%

\ifthenelse{\isodd{\value{page}}}{%
\makeatletter\normalmarginpar%
\marginpar{\tiny {#1}}%
}{%
\makeatletter\reversemarginpar%
\marginpar{\tiny {#1}}%
}}%

\usepackage{url}
\usepackage{color}
\usepackage{bm}
\usepackage{ulem}
\usepackage{comment}
\usepackage{siunitx}
\usepackage[]{defjs3}
\usepackage{siunitx}
\usepackage{lipsum}
\usepackage{acronym}
\usepackage{amsfonts}
\usepackage{graphicx}
\usepackage{epstopdf}
\usepackage{algorithmic}
\usepackage{multicol}
\usepackage{multirow} 
\ifpdf
  \DeclareGraphicsExtensions{.eps,.pdf,.png,.jpg}
\else
  \DeclareGraphicsExtensions{.eps}
\fi

\newcommand{\numgru}{\num[group-minimum-digits=3]}
\acrodef{DOFs}{degrees of freedom}
\acrodef{ROI}{region of interest}
\acrodef{GMRES}{Generalized Minimal Residual}
\acrodef{GDSW}{Generalized Dryja-Smith-Widlund}
\acrodef{PDEs}{partial differential equations}
\acrodef{PDE}{partial differential equation}
\acrodef{DDMs}{domain decomposition methods}
\acrodef{CTW}{Controlled Tensile Weldability}
\acrodef{LBW}{laser beam welding}
\acrodef{HPC}{High-Performance Computing}
\acrodef{FOV}{Field of view}

\newsiamremark{remark}{Remark}
\newsiamremark{hypothesis}{Hypothesis}
\crefname{hypothesis}{Hypothesis}{Hypotheses}
\newsiamthm{claim}{Claim}

\headers{Large-scale thermo-mechanical simulation of laser beam welding}{Large-scale thermo-mechanical simulation of laser beam welding}
\title{Large-scale thermo-mechanical simulation of laser beam welding using high-performance computing: \\
A qualitative reproduction of experimental results}

\author{Tommaso Bevilacqua\thanks{Department of Mathematics and Computer Science, University of Cologne, Germany (\email{tommaso.bevilacqua@uni-koeln.de, axel.klawonn@uni-koeln.de, martin.lanser@uni-koeln.de}).}
\and Andrey Gumenyuk\thanks{Bundesanstalt f\"ur Materialforschung und -pr\"ufung (BAM), Berlin, Germany (\email{andrey.gumenyuk@bam.de, niloufar.habibi@bam.de, michael.rethmeier@bam.de})}
\and Niloufar Habibi\footnotemark[2]
\and Philipp Hartwig\thanks{Institute of Mechanics, Department of Civil Engineering, University Duisburg-Essen, Germany (\email{philipp.hartwig.ph97@uni-due.de,j.schroeder@uni-due.de})}
\and Axel Klawonn\footnotemark[1] \thanks{Center for Data and Simulation Science, University of Cologne, Germany}
\and Martin Lanser\footnotemark[1] \footnotemark[4]
\and Michael Rethmeier\footnotemark[2]
\and Lisa Scheunemann\thanks{Chair of Applied Mechanics, Department of Mechanical and Process Engineering, RPTU Kaiserslautern-Landau, Germany (\email{lisa.scheunemann@rptu.de})}
\and J\"org Schr\"oder\footnotemark[3]}

\usepackage{amsopn}

\ifpdf
\hypersetup{
  pdftitle={Large-scale thermo-mechanical simulation of laser beam welding using high-performance computing: A qualitative reproduction of experimental results},
  pdfauthor={T. Bevilacqua, A. Gumenyuk, N. Habibi, P. Hartwig, A. Klawonn, M. Lanser, M. Rethmeier, L. Scheunemann and J. Schr\"oder}
}
\fi

\begin{document}

\maketitle
% REQUIRED
\begin{abstract}
Laser beam welding is a non-contact joining technique that has gained significant importance in the course of the increasing degree of automation in industrial manufacturing. This process has established itself as a suitable joining tool for metallic materials due to its non-contact processing, short cycle times, and small heat-affected zones. One potential problem, however, is the formation of solidification cracks, which particularly affects alloys with a pronounced melting range. Since solidification cracking is influenced by both temperature and strain rate, precise measurement technologies are of crucial importance. For this purpose, as an experimental setup, a Controlled Tensile Weldability (CTW) test combined with a local deformation measurement technique is used.

The aim of the present work is the development of computational methods and software tools to numerically simulate the CTW. The numerical results are compared with those obtained from the experimental CTW. In this study, an austenitic stainless steel sheet is selected. A thermo-elastoplastic material behavior with temperature-dependent material parameters is assumed. The time-dependent problem is first discretized in time and then the resulting nonlinear problem is linearized with Newton's method. For the discretization in space, finite elements are used. In order to obtain a sufficiently accurate solution, a large number of finite elements has to be used. In each Newton step, this yields a large linear system of equations that has to be solved. Therefore, a highly parallel scalable solver framework, based on the software library PETSc, was used to solve this computationally challenging problem on a high-performance computing architecture. Finally, the experimental results and the numerical simulations are compared, showing to be qualitatively in good agreement.
\end{abstract}

% REQUIRED

% REQUIRED

\section{Introduction}
%\begin{itemize}
%\item Something on laser beam welding and solidification cracks forming in general. 
%\item something on the specific experiment we consider here and why this setup is chosen / useful
%\item something on why simulating such experiments with FE and thermo-elastoplasticity is beneficial. A motivation for the paper essentially
%\item and then we might proceed as follows...
%\end{itemize}
Over the past decade, the \ac{LBW} process has undergone significant advancements and has established itself as an efficient and economical tool in the industry. This process offers a highly concentrated energy input, minimizing thermal strain in components while enabling rapid processing with deep penetration—ideal for producing durable, high-performance parts. Key applications can be found in high-precision industries such as aerospace, energy technology, and plant and apparatus engineering. Temperature- and corrosion-resistant high-alloy steels are commonly used in these safety-critical areas. However, during welding, these materials are prone to solidification cracking — a complex phenomenon influenced by the interplay of various thermal, metallurgical, and mechanical factors \cite{agarwal2019study,lippold1994solidification, norouzian2023review}. This type of cracking typically occurs during the final stages of solidification in the mushy zone - the region where both solid and liquid phases coexist -, when tensile stresses exceed the ability of the remaining liquid to backfill gaps. Additionally, high solid fractions and an insufficient liquid flow further intensify cracking. 

Therefore, an accurate determination of the ductility limits of a material and the influence of process parameters is essential for the optimization of the LBW process. This requires experimental methods to evaluate ductility in the mushy zone, that is within the brittle temperature range (BTR) \cite{prokhorov1962resistance}, which spans from the solidus to the liquidus temperature. Since solidification cracking is influenced by both temperature and strain rate, precise measurement technologies are crucial along with a testing system, which is capable of applying varied strain conditions during welding. To this end, \ac{CTW} \cite{quiroz2012investigation}, combined with high-resolution temperature and local strain measurement technologies, is an effective method, allowing, in particular, for the application of a wide range of global strains and strain rates during welding. Additionally, integrating computer-aided engineering in parallel with experimental tests provides deeper insights into the process and enhances the accuracy of ductility limit determinations, which can subsequently be utilized in the design of complex manufacturing processes \cite{Bakir2020}.

Here, our goal is to set up and present efficient finite element simulations of the \ac{CTW} test to get new insights using highly resolved simulations, exploiting the potential of modern supercomputers. First numerical methods for the description of thermo-mechanical problems for welding were presented in \cite{Fri:1975:tao} and \cite{HibMar:1973:ant} and nowadays, finite element simulations of \ac{LBW} processes can, in general, be performed in various software environments, each of them with several advantages and disadvantages. On the one hand, simulations in commercial finite element software tools like ANSYS~\cite{bakir2024determination,ranjbarnodeh2016finite} can be set up quite fast and parameters such as problem geometries, boundary conditions, finite element types, or material models can be varied and tested easily. Also, important features such as a non-uniform mesh with local refinements close to the weld seam can be implemented without much effort. On the other hand, customized finite element simulation tools can be used which are able to perform large-scale computations on modern supercomputers. They can deliver high-definition solutions in a fraction of the computing time and thus new insights in the whole process of \ac{LBW}. Unfortunately, the overhead of mirroring the exact experiment in such a software framework is usually high. As a consequence, we decided to set up a hybrid software pipeline exploiting the advantages of both worlds: the easy-to-use commercial software ANSYS and the high-performance software FE2TI~\cite{klawonn2020computational}, which is based on the very efficient PETSc (Portable, Extensible Toolkit for Scientific Computing) library~\cite{balay2019petsc}. More precisely, we exploit ANSYS to perform a finite element simulation on the complete geometry mirroring, as close as possible, the experimental set up and extract the deformation and temperature data over time on the boundary of a small strip around the weld seam. Afterwards, we load these boundary data into FE2TI, interpolate them to a finer mesh and perform high-performance computational simulations on modern supercomputers, specifically of the area around the weld seam - the actual zone of interest. These simulations deliver highly resolved and accurate results within hours which is impossible to obtain with tools like ANSYS directly. Let us remark that in ANSYS a sequential approach is used where first only the temperature field over time is computed and afterwards a temperature-driven elastoplasticity problem is solved. In FE2TI, a fully coupled thermo-elastoplasticity problem is solved leading to a more realistic coupling between temperature and deformation. The coupled thermo-elastoplastic formulation is based on \cite{SimMie:1992:act}, in which finite strains were taken into account, and adapted to small strains. The plastic behavior is modeled by von Mises plasticity, see \cite{Mis:1913:nvd}, and takes into account multilinear hardening. The computations are based on trilinear hexahedral elements for both displacement and temperature fields, see \cite{SimTayPis:1985:vap}.

Overall, the goal of this work is to represent an initial qualitative comparison between experiments and numerical simulations, as well as to highlight their challenges, issues, and possible solutions. So far, we noticed qualitative similarities between these two approaches; although a quantitative comparison is not yet possible, this study still provides valuable insights. In particular, we will see how numerical simulations can be a complement to the experiments, which are limited to analyzing the surface of the plate and cannot investigate what happens inside, where crack initiation appears to take place.

Regarding the \ac{HPC} framework in FE2TI, we use parallel domain decomposition solvers, which are designed and optimized to handle large thermo-elastoplasticity finite element problems~\cite{bevilacqua2025highly,bevilacqua2024monolithic}. More precisely, we use overlapping Schwarz \ac{DDMs}~\cite{toselli2006domain} with variations of GDSW (Generalized Dryja-Smith-Widlund) coarse spaces which already showed to be robust for coupled multi-physics applications in the past; see~\cite{dohrmann2008family,dohrmann2008domain,heinlein2019monolithic} for details. Several additional features, such as, recycling strategies in the solver setup, adaptive time stepping, the optimal choice of sparse direct solvers, and optimized stopping criteria are exploited to further decrease the simulation time as well as the energy consumption of the simulation and to increase its computational efficiency~\cite{bevilacqua2025highly}.

%\todo[inline]{ML: If you ask me, we can delete the ''remainder'' part. I usually find it useless.}
The remainder of the article is organized as follows. In~\cref{sec:experiment}, we provide a detailed description of the experimental setup. This is followed by a description of the finite element and material models used in all computations and the exact description of the boundary value problems in ANSYS and FE2TI in Sections~\ref{sec:FE-model} and~\ref{sec:mat_parameter}. Convergence studies, scalability results, and a comparison of the simulation data obtained from the FE2TI simulations with the experimental results are then considered in~\cref{sec:hpc}. There, also a detailed description of the solver setup of the chosen \ac{DDMs} in FE2TI is given. We finish with some concluding remarks in~\cref{sec:conclusion}.

\section{Experimental setup} 
\label{sec:experiment}
% Introducing sentences for each section should be deleted (Meeting 5.3.2025)
%{\TB In this section we describe the material used into our simulation and further introduce the \ac{CTW} experiment and the relative optical measurement technique utilized to analyze the simulation.}

\subsection{Material}
For our study, we select an austenitic stainless steel sheet of AISI 304 (1.4301) grade with a thickness of 1 mm and a chemical composition as given in \cref{tab:chem_compo}. Since the presence of the primary ferrite phase increases the sensitivity to hot cracking, the ferrite content is measured through metallographic analysis, the results of which indicate a ferrite content of 1.1\% for the base metal and 18\% for the fully penetrated weld seam. Generally, austenitic steels with  $\delta$-ferrite content ranging from 5\% to 15\% are considered to exhibit good resistance to hot cracking during welding \cite{lienert2003improved}. The sensitivity of this material to solidification cracking is therefore evaluated using the \ac{CTW} test, as detailed in the following section.
\renewcommand{\arraystretch}{1.0}
\begin{table}[H]
\centering
\caption{Chemical composition of the studied  CrMnNi stainless steel $(1.4301)$.}
    \begin{tabular}{|l|l|l|l|l|l|l|l|l|l|}
    \hline
    Element: & C    & Cr     & Mn   & Ni    & Si  & P     & S    & N    & Fe   \\ 
    \hline
    wt$\%$: & $0{.}02$ & $19{.}09$ & $1{.}6$ & $8{.}06$ & $0{.}41$ & $0{.}028$ & $<0{.}002$ & $0{.}095$ & bal. \\ 
\hline
\end{tabular}
\label{tab:chem_compo}
\end{table}

\subsection{CTW test procedure}\label{sec:ctw_test}
\begin{figure}[]
    \centering
    \begin{minipage}{0.49\textwidth}
        \centering
        \includegraphics[width=\textwidth]{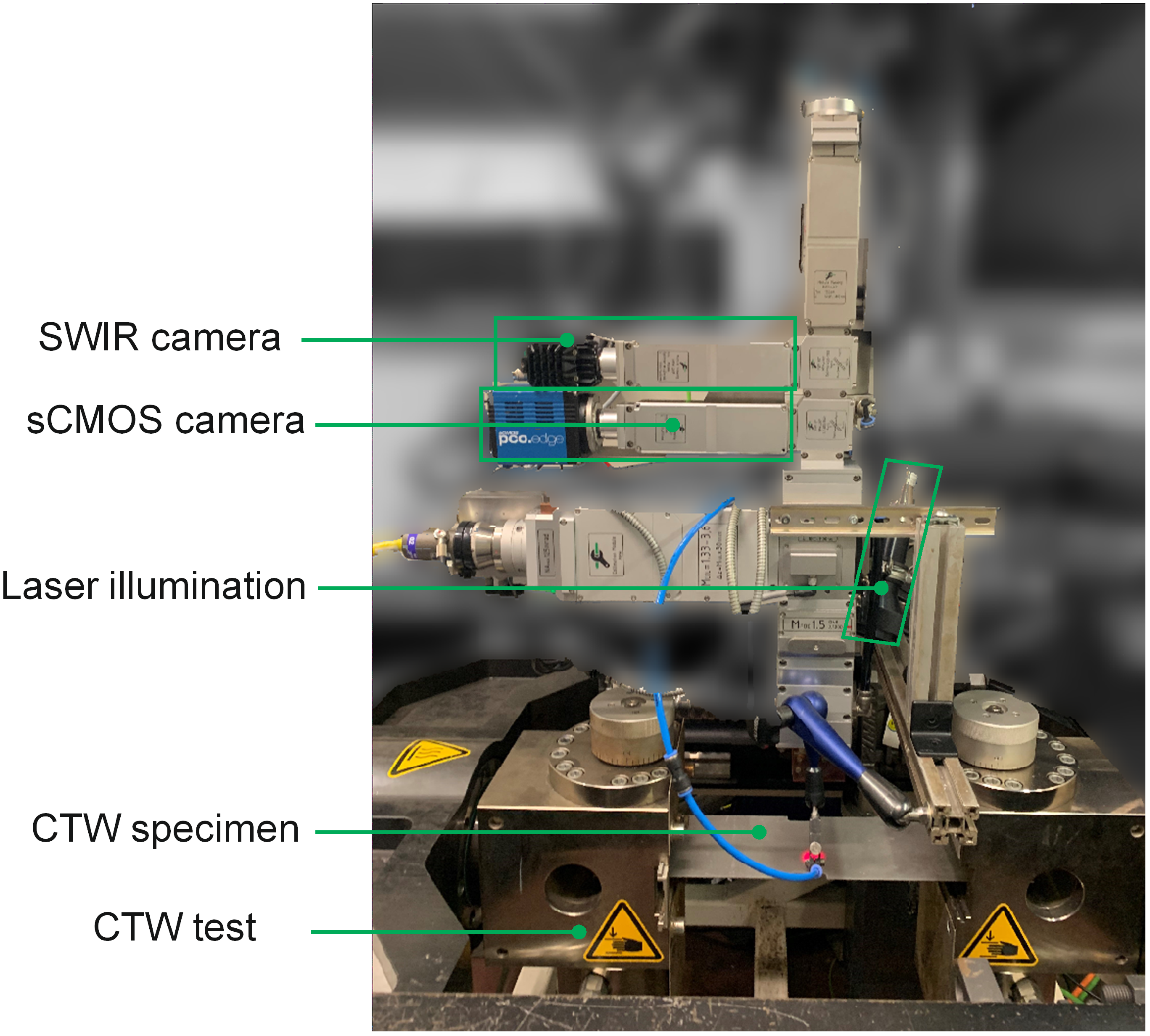}
    \end{minipage}
    \hfill
    \begin{minipage}{0.49\textwidth}
        \centering
        \includegraphics[width=\textwidth]{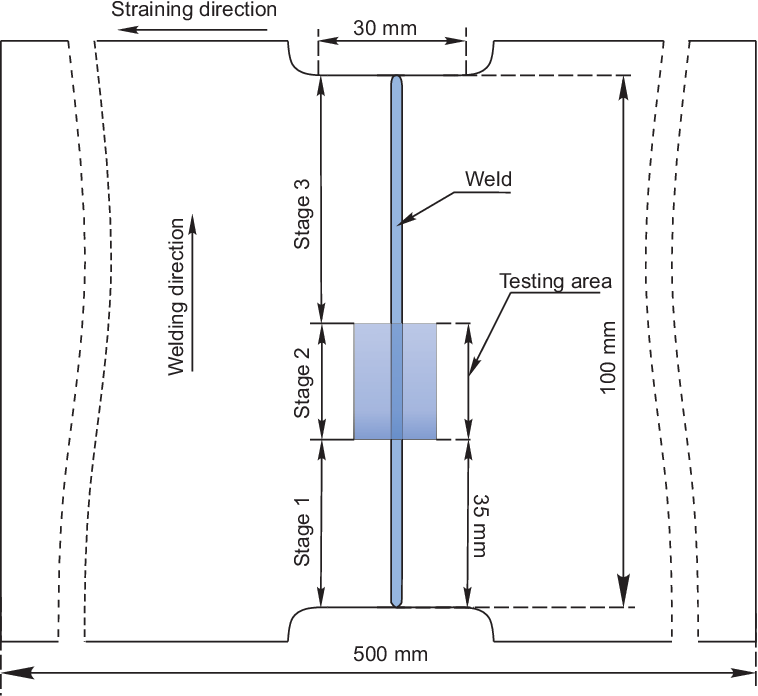}
    \end{minipage}
    \caption{(Left) Instrumented \ac{CTW}-test. (Right) Welding experiments under \ac{CTW}-Test conditions-specimens arrangement.}
    \label{fig:ctw_combined}
\end{figure}

The \ac{CTW} test is a hot cracking evaluation method developed at the Bundesanstalt für Materialforschung und -prüfung (BAM) in Berlin \cite{Coniglio2008AluminumAW}. It employs a horizontal tensile testing apparatus with a $50\,\text{kN}$ load capacity to impose a controlled planar tensile strain perpendicular to the welding direction during welding. This method involves systematically varying the global strain rate (crosshead speed) to influence the local strain rate near the mushy zone. By evaluating crack and no-crack conditions, the test aims to determine the critical global strain and strain rate required for crack formation. In this study, as shown in \cref{fig:ctw_combined} Left, the \ac{CTW}-setup is equipped with a high-resolution sCMOS camera (see \cref{sec:opt_meas_tech}) and two-color infrared thermograph (see \cite{gumenyuk2024two} for further details) to capture local deformation and temperature, respectively.
%\footnote{Detailed information in \cite{gumenyuk2024two}. However, the results are not used in the present work.}

The sheet is cut into $120\,\text{mm}\times500\,\text{mm}$ specimens (\cref{fig:ctw_combined} Right), and the manufactured notches at the center ensure the applied strain localizes in this region, that is a designed gauge area of $100\,\text{mm}\times 30\,\text{mm}$. 
The test consists of three stages. In the first one, the welding process is stabilized, and a consistent temperature field is established, while no external strain is applied. The characteristics of the laser and the welding parameters are listed in \cref{table:Table2}. The second stage begins after $35\,\text{mm}$ of welding, at which point external loading is applied at a constant $6\,\%/\mathrm{s}$ global strain. Once the maximum global strain is reached, the external loading ceases, and the welding process continues until a total weld length of $90\,\text{mm}$, the third stage. For this study, the maximum global strain varies from $0.02$ up to $0.035$ with the step size of $0.005$.

\begin{table}[H]
\centering
	\caption{Yb:YAG disc laser characteristics and the welding parameters}
	\label{table:Table2}
	\begin{tabular}{@{}crcr@{}}
		\hline
		Max. nominal power in $\mathrm{kW}$ & 16   \vline &Welding speed in $\mathrm{m/min}$ & 1 \\ \hline 
		Wavelength in $\mathrm{nm}$ & 1030  \vline  &Laser power in $\mathrm{kW}$ & 1  \\ \hline 
		Fiber core diameter in $\mathrm{\mu m}$ & 200  \vline &Ar flow rate in $\mathrm{1/min}$ & 20\\ \hline 
		\begin{tabular}{c}		
		Beam parameter product (BPP) \\
		in $\mathrm{mm \times mrad}$		
		\end{tabular}
		 & 8 \vline & Focal position in $\mathrm{mm}$ & 5 \\ \hline 
		Focal distance in $\mathrm{mm}$ & 300  \vline  &Focus diameter in $\mathrm{\mu m}$ &450  \\ \hline 
	\end{tabular}
\end{table}

\subsection{Optical measurement technique}\label{sec:opt_meas_tech}

\begin{figure}[]
\centering
\includegraphics[width=0.8\textwidth]{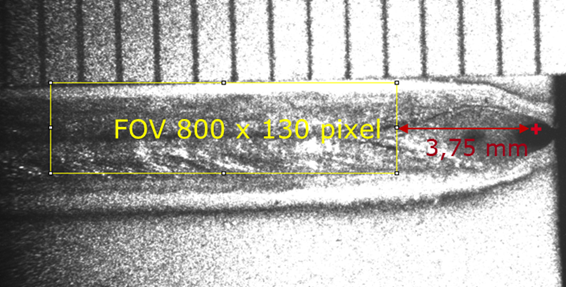}
\caption{\ac{FOV} of the camera image}\label{fig:FOV}
\end{figure}

The welding process is recorded using a pco.edge 5.5 sCMOS camera, which is installed co-axially to the laser path, see \cref{fig:ctw_combined} Left. To ensure a valid recording of the welding process and the re-solidified material, an external $100\,\text{W}$ diode laser with a central wavelength of $808\,\text{nm}$ is used for illumination. In this setup, the illuminated spot size corresponds roughly to an ellipse with dimensions of $30\,\text{mm}\times60\,\text{mm}$. A narrow bandpass interference filter of the same wavelength and a bandwidth of $20\,\text{nm}$ is placed in front of the camera, allowing only the wavelengths from the illumination to pass through and suppressing process emissions in all other spectral ranges during recording. The recording is conducted at $778\,\text{fps}$, with images stored in TIFF format at a resolution of $800\,\text{pixels}\times130\,\text{pixels}$. The captured area covers a region of $9\,\text{mm}\times1.46\,\text{mm}$, positioned $3.75\,\text{mm}$ behind the weld pool along the laser beam axis, as illustrated in \cref{fig:FOV}.

The local straining condition is calculated at the immediate vicinity of the solidification front using a novel optical technique, based on the Optical Flow (OF) along with Lukas-Kanade (LK) algorithm \cite{bakir2019development}. Because the laser head, where the camera is inserted, typically moves while the workpieces remain stationary, only the rear part of the weld pool is evaluated. Furthermore, since solidification cracking occurs in the mushy zone, focusing on this area as the \ac{ROI} is reasonable. This means that the \ac{ROI} changes at each step, thus all the related calculations must be updated correspondingly. Therefore, the challenge in this task is tracking of the pixels in the critical temperature range to obtain the displacement, which will be used later to calculate the strain. For the \ac{ROI}, an area of $430\,\text{pixels}\times130\,\text{pixels}\ (4.8\,\text{mm} \times 1.46\,\text{mm}$) is considered.

The \ac{ROI} can be imagined as a mass balance for the steady-flow process that has one inlet and one outlet. The material pixels, during solidification, will enter the evaluation region from the inlet boundary and be tracked until they exit from the opposite boundary or the outlet. The new pixels, which exceed the edge of the \ac{ROI} due to the movement of the material, take the displacement values of 0; then, they can be considered in the displacement computation. After estimating the displacement between two sequential frames, the main displacement field can be warped according to the newly computed displacement and then added to the temporary displacement field. The current algorithm is updated to provide a more robust computation of displacements by means of OF algorithms accounting for averaging over several frames. Based on the estimated displacements, the Green-Lagrangian strain can be computed for each frame.

%--- Section ---%
\section{Finite Element  Modeling}\label{sec:basic_set_eq}\label{sec:FE-model}
%In this section we describe the material modeling, the finite element approach and the heat soure modeling used into the ANSYS and FE2TI simulations.
%\todo[inline]{TODO: A nice couple of sentences to justify the fact that we are using different but sufficiently similar FEM discretization and plastic models}

\subsection{Material model in ANSYS}
In ANSYS, the material behavior is modeled using a multilinear isotropic hardening plasticity model based on the von Mises yield criterion (see \cref{fig:yield_curves}). Initially, the material exhibits linear elastic behavior, defined by a temperature-dependent Young modulus and Poisson ratio. Upon reaching the yield point, the material transitions into plastic deformation, described by piecewise linear hardening curves for each temperature. This approach captures the nonlinear stress-strain response and accounts for the material's strength degradation at elevated temperatures, making it suitable for simulating ductile materials under complex loading and thermal conditions. 

\subsection{Fully coupled thermo-elastoplasticity at small strains (FE2TI)}\label{sec:sub_basic_eq}

%{\ML All large scale simulations are carried out using FE2TI~\cite{klawonn2020computational}, a highly scalable and PETSc-based~\cite{balay2019petsc} solver framework for single- or two-scale Finite Element simulations of solid mechanics problems; see~\cref{sec:hpc} for details. The phenomenological material models itself are implemented in the finite element library FEAP~\cite{Tay:2008:fea} and provided to FE2TI through an efficient and lightweight software interface connecting both software packages.}
To perform the large finite element simulations carried out on our FE2TI library, we equipped it with an efficient and lightweight software interface in order to take the phenomenological material models implemented in the finite element library FEAP~\cite{Tay:2008:fea}. This allows to combine the flexibility of FEAP in defining complex material models and our highly scalable, PETSc-based solver framework.\\
More precisely, for the simulation of \ac{LBW} processes and the \ac{CTW} tests in FE2TI, a material model was implemented in FEAP that mimics the thermo-elastoplastic material behavior, see~\cite{SimMie:1992:act}, at small strains. Additionally, a hexahedral finite element with $\bar{\bB}$-formulation, see~\cite{SimTayPis:1985:vap}, was added. 
For the thermo-elastoplastic model, an additive decomposition of the total strains into an elastic part $\Bve^e$, a plastic part $\Bve^p$, and a thermal part $\Bve^t$ is assumed, such that $\Bve=\Bve^e+\Bve^p+\Bve^t$. Thus, the applied free energy function $\psi$ can be formulated, which is dependent on the thermo-elastic strains $\Bve^{te}=\Bve^e+\Bve^t$, the temperature $\theta$, and the strain-like internal variable  $\alpha$, i.e.

%FEAP offers flexibility in defining complex material models (e.g., nonlinear, anisotropic, and user-defined materials) and supports various FEM formulations, including solid mechanics, thermal analysis, and multiphysics simulations.

\begin{equation}
    \begin{split}
    \psi\left(\Bve^{te},\theta,\alpha\right)&=\frac{1}{2}\kappa\left(\tr\Bve^{te}\right)^2+\mu\Vert\dev\Bve^{te}\Vert^2-c_\rho\left(\theta\ln{\frac{\theta}{\theta_0}}-\theta+\theta_0\right)\\
    &-3\alpha_T\kappa\left(\theta-\theta_0\right)\tr\left(\Bve^{te}\right)+f(\alpha)
    \end{split}
\end{equation}
comp. \cite{SimMie:1992:act}. Here, $\kappa$ denotes the bulk modulus, $\mu$ is the shear modulus, $c_\rho$ is the heat capacity, $\theta_0$ is the initial temperature, $\alpha_T$ is the heat expansion coefficient, and $f$ is the hardening function. In addition to the free energy function the balance of momentum  
\begin{equation}
\mathrm{div}\boldsymbol{\sigma}+\bbf=\boldsymbol{0}\ \text{with}\ \boldsymbol{\sigma}=\boldsymbol{\sigma}^T
\label{equ:bal_mom}
\end{equation}
and the balance of energy
\begin{equation}
\rho r-\div\bq+\rho\theta\frac{\partial^2\psi}{\partial\theta\partial\Bve}:\dot{\Bve}+\rho\theta\frac{\partial^2\psi}{\partial\theta^2}\dot{\theta}=0
    \label{equ:bal_energy}
\end{equation}
need to be considered. In this case, $\bbf$ describes the body forces and $\bq$ describes the heat flux vector. The application of the entropy inequality results in the constitutive equations $\Bsigma=\p{\Bve^{te}}\psi$, which describes the relationship between stresses and strains, and $\beta=\p{\alpha}\psi$, which describes the relationship between the strain-like internal variable $\alpha$ and the stress-like variable $\beta$. To describe the plastic behavior, a von Mises yield criterion

\begin{equation}
    \Phi=\Vert\dev\Bsigma\Vert-\sqrt{\frac{2}{3}}\left(y_0+\beta\right),
\end{equation}

see~\cite{SimHug:1998:cin}, is applied where $y_0$ denotes the yield strength, which is a function of the temperature and the accumulated plastic strains. \\
%\todo[inline]{Maybe we can be short here, since finite element discretization is well known? I make a suggestion:}
The \ac{PDEs} are time-dependent and nonlinear. We use an implicit Euler scheme to discretize the nonlinear system in time and Newton's method to linearize it in each time step. Let us remark that the resulting coupled linear systems have a two-by-two block structure where the first diagonal block represents the mechanical problem and the second one the thermal diffusive problem. The coupling is nonlinear and not symmetric. Besides usual boundary conditions of Dirichlet- and Neumann-type for displacements and temperature, we have to implement the heat entry of the laser in the system, which is described in the following section.

\subsection{Heat Source modeling}\label{sec:heat_source}
\begin{figure}[]
\centering
    \includegraphics[clip=true, trim={500 400 500 390},width=0.7\textwidth]{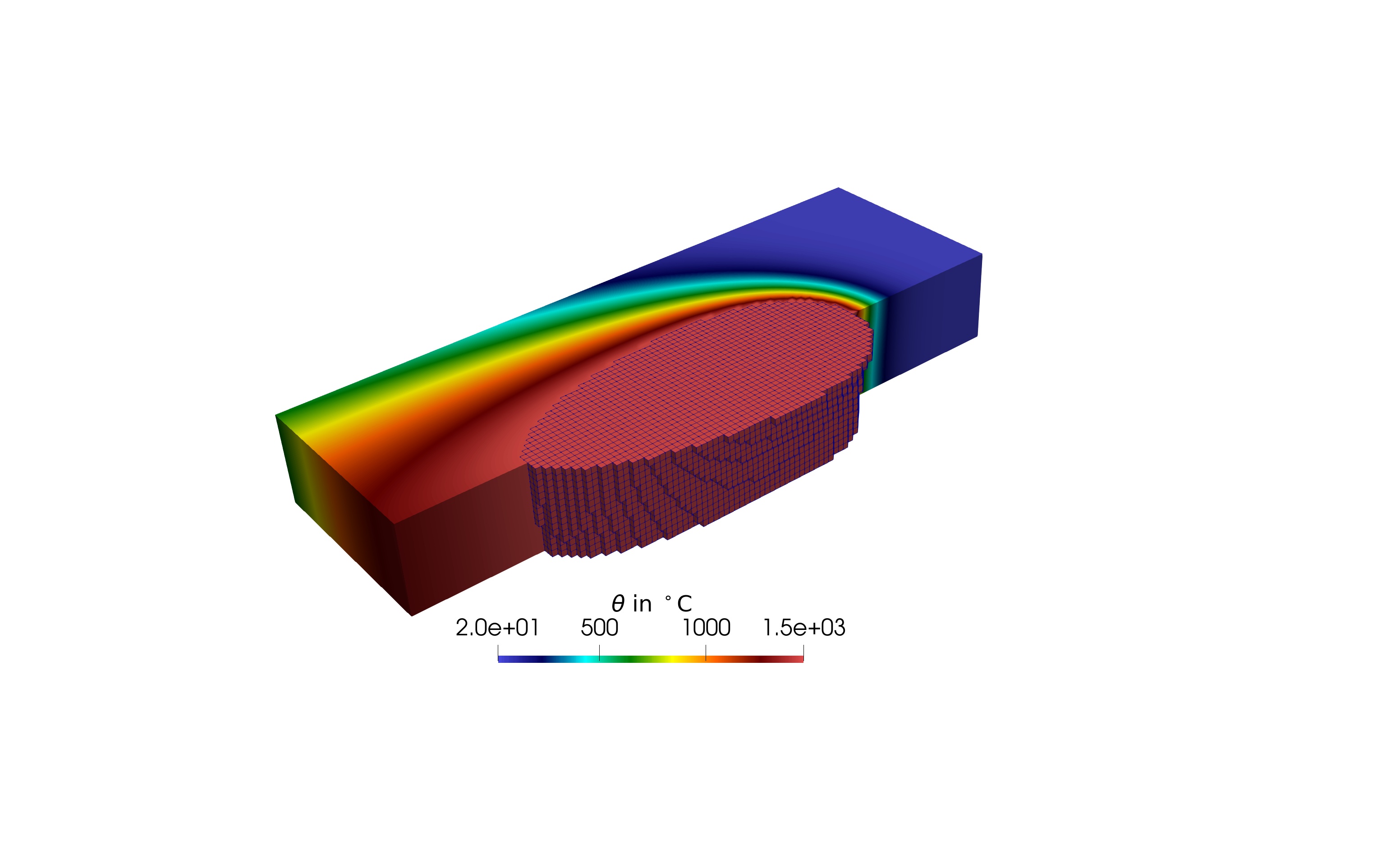}
\caption{Temperature field and discrete representation of the melting pool using hexahedral finite elements.}\label{fig:melt_pool_geom}
\end{figure}

The action of the laser is modeled as a constraint on a volume over the temperature field, that represents the molten region of the metal (melting pool); see \cref{fig:melt_pool_geom}. To achieve this, we utilize a surface model derived from experimental data to represent this region and we move this geometry along the welding direction during the time iterations at a velocity of $1.0 \, \mathrm{m}/\mathrm{min}$, see \cref{table:Table2}. The initial temperature of the metal is set uniformly to $\theta_0 = 20^\circ $C and the temperature within the melting pool is set to the liquidus temperature of the material of $\theta_\mathrm{liq} = 1460^\circ $C.
While ANSYS allows the initial temperature \ac{DOFs} to be instantly set to $\theta_\mathrm{liq}$, in the FE2TI framework, this must be enforced through a linear increment for stability reasons. A heating procedure (initialization phase) is therefore implemented to ensure that the temperature $\theta_\mathrm{liq}$ is reached within the melting pool after $0.1 \, s$. After this initialization phase, the movement of the melting pool starts and it is necessary to constantly update the temperature constraints, since in each time step some finite elements enter the melting pool and some leave it. In contrast to the initialization phase, when a temperature \ac{DOFs} enters into the moving melting pool is immediately set to $\theta_\mathrm{liq}$, which caused no problems in the simulations so far.

%--- Section ---%
\section{Boundary Value Problem  and Material Parameter}\label{sec:mat_parameter}

%In this section we provide a detailed description of the finite element setup of the considered \ac{CTW} test and also report all necessary temperature-dependent material parameters used in the thermo-elastoplastic material model.

\subsection{Boundary Value Problem} 
    
\begin{figure}[]
\centering
\includegraphics[width=1\textwidth]{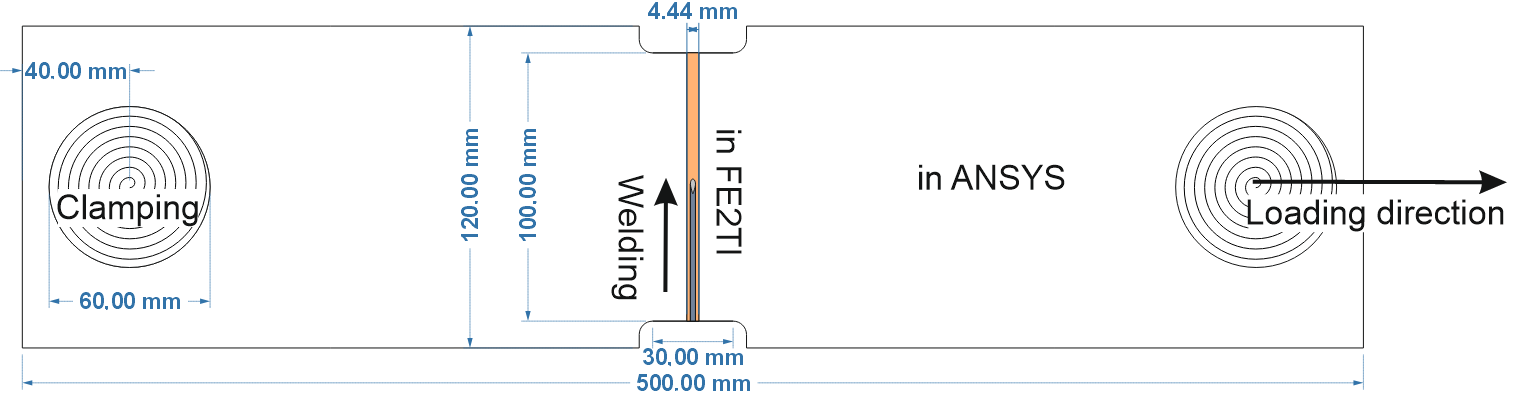}
\caption{Schematic illustration of \ac{CTW} test; the whole specimen was used for the ANSYS simulations, while the FE2TI simulations focus on the middle strip.}\label{fig:bvp_ansys}
\end{figure}

The boundary value problem at hand, which is solved in the simulations, is based on the experimental setup of the \ac{CTW} test. In the ANSYS simulations the entire system is used, while a section of the specimen is considered for the simulation with FE2TI. The entire specimen therefore has dimensions of $500\,\mathrm{mm}$ by $120\,\mathrm{mm}$ and a thickness of $1\,\mathrm{mm}$. In the center of the system there is a cut-out over a length of $30\;\mathrm{mm}$, which is tapered by $10\,\mathrm{mm}$ on both sides to $100\,\mathrm{mm}$. This is intended to concentrate the arising strains in the area of the cut-out; see \cref{fig:bvp_ansys}.

\begin{figure}[]
\centering
\includegraphics[width=0.7\textwidth]{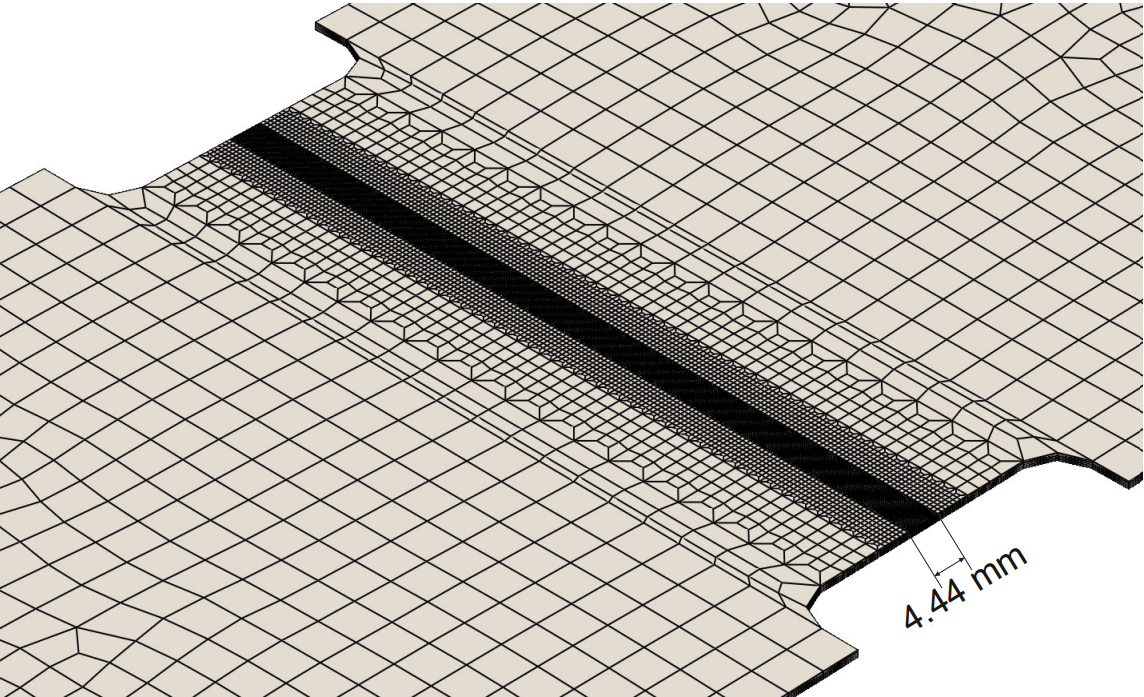}
\caption{The mesh pattern used in the ANSYS simulation; the results on the boundary of the marked $4.44\,\mathrm{mm}$ strip are written out during the simulation and employed in FE2TI simulations as boundary conditions.}\label{fig:bvp_ansys_mesh}
\end{figure}

For the simulation in ANSYS, a thermal analysis is first performed to determine the temperature history. Then, the temperature results are applied incrementally to the mechanical model to calculate stresses and strains. A step size of 0.001\,\text{s} is defined for both processes. The system is meshed using  88\,780 elements with a finer size of $0.185\,\text{mm}\times0.185\,\text{mm}\times0.2\,\text{mm}$ at the critical region. The resulting finite element mesh is shown in \cref{fig:bvp_ansys_mesh}. For this regard, the element types of SOLID70 (3-D Thermal Solid) and SOLID186 (3-D 20-Node Structural Solid) for the thermal and mechanical analyses are used, respectively.

To simplify the computation, the laser itself is not simulated, but the heat source model described in Section \ref{sec:heat_source} is used. More precisely, we recall that the isotherm of the liquidus temperature $(1460^\circ\mathrm{C})$ is predefined and moved forward along the welding direction identically to the real welding process. A constant temperature therefore prevails within the isotherm, while the temperature outside the isotherm is freely adjustable. The mechanical boundary conditions are identical to those of the experiment. For this purpose, two circular areas with a diameter of $60\,\mathrm{mm}$ are defined at opposite ends in which the boundary conditions are applied. The two areas are centered along the welding direction and have a distance between center and the outer edge of the specimen of $40\,\mathrm{mm}$. In one area, all displacements are prevented so that a clamping is present, while in the second area the displacements in the $x$- and $z$-directions are prevented and a prescribed displacement is applied for the $y$-direction. For the applied displacements, a maximum strain of $2.5\,\%$ with a strain rate of $6\,\%/\mathrm{s}$ is assumed. However, the strains are not applied right from the start, but increase with the strain rate from the point at which the laser has traveled a distance of $35\,\mathrm{mm}$. Once the maximum strains have been reached, they are not increased any further but remain constant until the end of the welding process. This mimices the behavior of the experiment as described in \cref{sec:ctw_test}.

\begin{figure}[]
\centering
\includegraphics[clip=true, trim={0 100 0 460},width=1\textwidth]{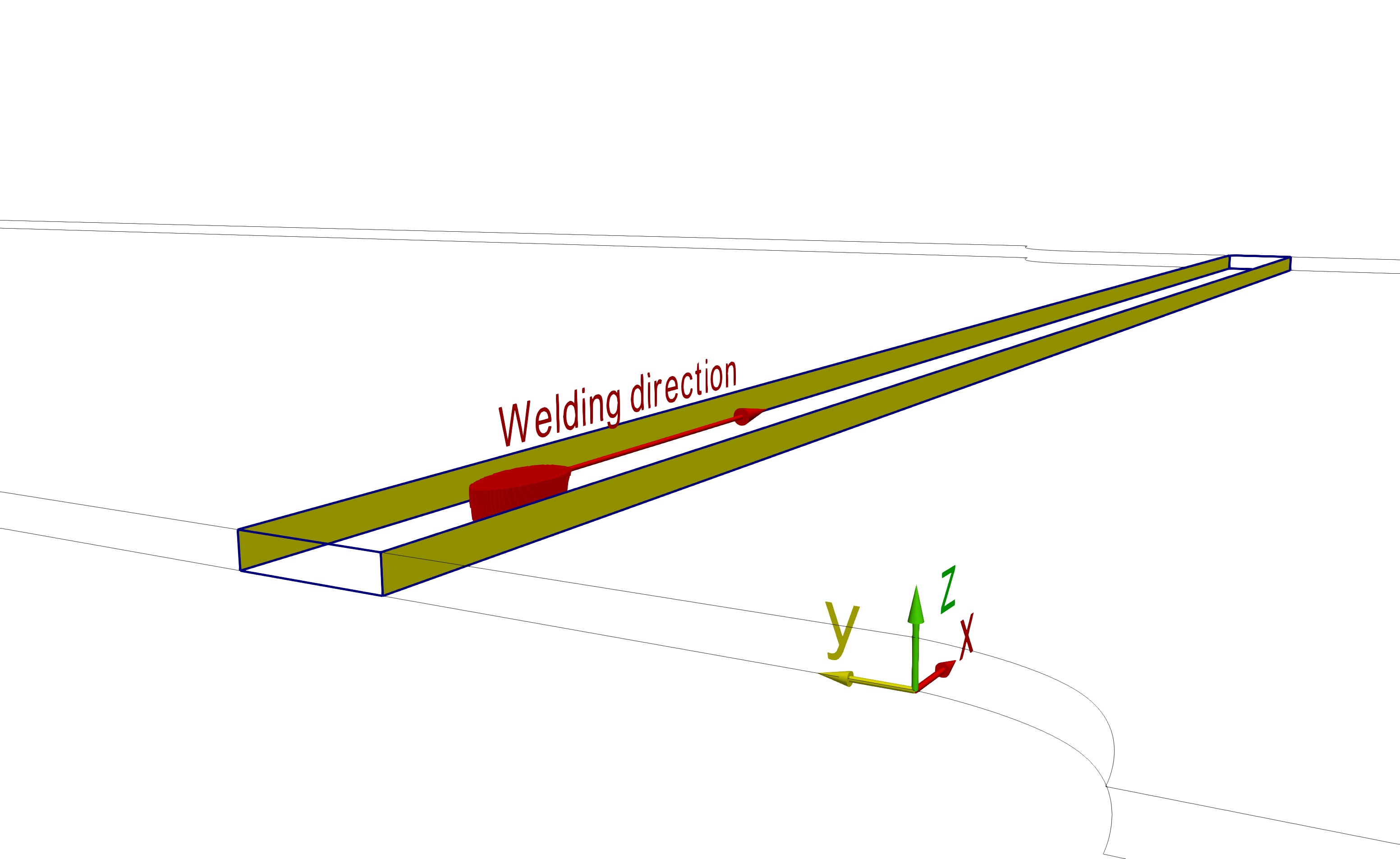}
\caption{Schematic illustration of the domain for the FE2TI simulation. In yellow the surfaces where the boundary conditions from ANSYS are applied.}\label{fig:bvp_fe2ti}
\end{figure}

As already mentioned, the simulation in FE2TI is carried out on a sub-geometry of the simulation with ANSYS which results in a cuboid with a length of $100\,\mathrm{mm}$, a width of $4{.}44\,\mathrm{mm}$ and a height of $1\,\mathrm{mm}$; see \cref{fig:bvp_ansys}. Trilinear structured hexahedral elements are used to create a highly resolved mesh with up to $4.4$ millions of elements. As the simulation with FE2TI only solves a subproblem, the mechanical and thermal boundary conditions are extracted from the simulation with ANSYS at the corresponding edges of the associated sub-geometry. These boundary conditions are then applied as Dirichlet boundary conditions on the lateral surfaces along the welding direction (see \cref{fig:bvp_fe2ti}). To achieve more accurate results, smaller time steps and finer meshes are used in FE2TI simulations. Therefore, a linear interpolation over time and a bilinear interpolation over space is used to enforce the boundary conditions on the finer discretization. Finally, the thermal constraints to represent the laser action are enforced as previously described in \cref{sec:heat_source}.

\subsection{Material Parameters}
\begin{figure}[]
\centering
\def\svgwidth{12cm}
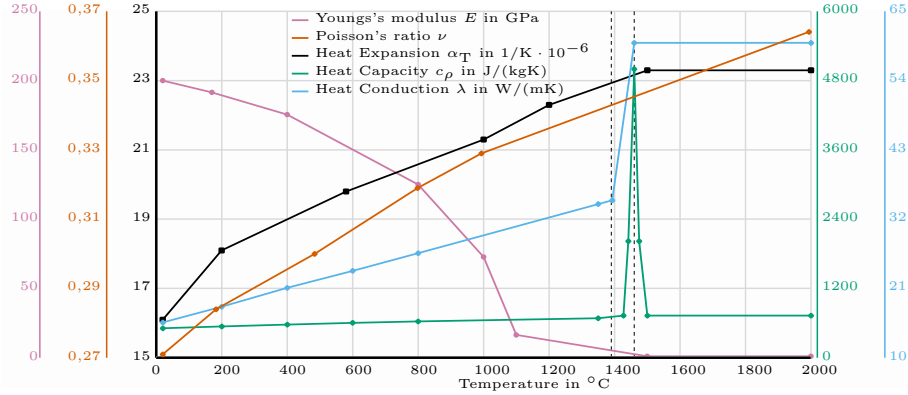
\caption{Material parameters for the used austenitic chrome-nickel steel 1.4301. The data originate from the software \textit{Sysweld} \cite{ESI:2009:Sys}. The first dashed vertical line visualizes the solidus temperature $\theta_\mathrm{sol}=1390^\circ\mathrm{C}$ and the second dashed line visualizes the liquidus temperature $\theta_\mathrm{liq}=1460^\circ\mathrm{C}$.}
\label{fig:material_parameter}
\end{figure}

% <<<<It is already mentioned in experimental section above >>>>>The material used in this paper is an austenitic chromium-nickel steel AISI 304.%with the material number 1.4301, also known as X5CrNi18-10. The chemical composition of the material used in the experiments can be found in Table \ref{tab:chem_compo}. 
The material parameters required for the simulation are the Young modulus $E$, the Poisson ratio $\nu$, the coefficient of thermal expansion $\alpha_T$, the specific heat capacity $c_\rho$, the thermal conductivity $\lambda$ and the density $\rho$. All these values, except for the density, are evaluated for a temperature range from $0^\circ\mathrm{C}$  to $2000^\circ\mathrm{C}$ with the help of the software \textit{Sysweld} \cite{ESI:2009:Sys}
%(https://www.esi-group.com/products/sysweld, access on {\AG06:03:2025})
and are shown in \cref{fig:material_parameter}.

% \subsection{Material Parameters}\label{sec:mat_para}

Since no density data is available in Sysweld, a density of $7919\,\mathrm{kg}/\mathrm{m}^3$ was taken from~\cite{Ric:2010:dpe} and it is assumed to be constant over temperature. Since the data is only available for some specific temperatures, a piecewise linear interpolation between two known values is used in order to determine the material parameters for any temperature value. We obtain

\begin{equation}
\theta_k < \theta \leq \theta_{k+1}\ \Rightarrow\  p =p_k+\dfrac{p_{k+1}-p_k}{\theta_{k+1}-\theta_k}\cdot\left(\theta-\theta_k\right)\ \text{with}\ p\in\{E,\nu,\alpha_\textrm{T}, c_\rho, \lambda\},
\label{equ:interpol}
\end{equation}
where $\theta_k$ and $\theta_{k+1}$ reflect the temperatures of a respective data point and $\theta$ the intermediate temperature at which the material parameters are evaluated. Accordingly, $p$ stands for the material parameter to be determined and $p_k$ and $p_{k+1}$ for the material parameters of the associated temperatures.

\begin{figure}[]
\centering
\def\svgwidth{12cm}
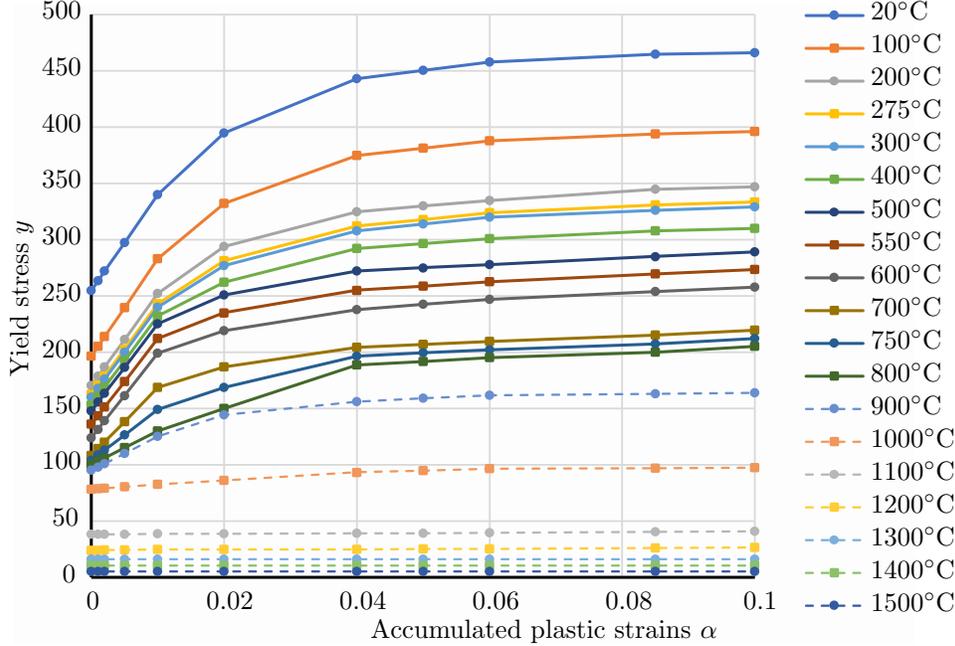
\caption{Yield curves for the used austenitic chrome-nickel steel 1.4301. The~data originate from the software \textit{Sysweld}.}
\label{fig:yield_curves}
\end{figure}

% \subsection{Yield Strength}\label{sec:yield_str}
Since a thermo-elastoplastic material behavior is assumed, values for the yield stress $y_0$ and the linear hardening parameter $h$ are needed. Both values depend on two variables, the temperature $\theta$ and the accumulated plastic strains $\alpha$. The available yield curves are shown in \cref{fig:yield_curves}.
As a result of the dependency on two variables, a two-step interpolation scheme is necessary. If the value of the yield stress and the linear strain hardening parameter has to be evaluated for a temperature $\theta_n$ and an accumulated plastic strain $\alpha_n$, the first step is to interpolate with respect to the temperature. Thus, $\theta_l<\theta_n<\theta_{l+1}$ applies and interpolation takes place at the points $\alpha_k$ and $\alpha_{k+1}$. It holds that

\begin{equation}
\begin{split}
y\left(\alpha_k,\theta_n\right)=y_0\left(\alpha_k,\theta_l\right)&+\frac{y_0\left(\alpha_k,\theta_{l+1}\right)-y_0\left(\alpha_k,\theta_l\right)}{\theta_{l+1}-\theta_l}\left(\theta_n-\theta_{l}\right) \\
y\left(\alpha_{k+1},\theta_n\right)=y_0\left(\alpha_{k+1},\theta_l\right)&+\frac{y_0\left(\alpha_{k+1},\theta_{l+1}\right)-y_0\left(\alpha_{k+1},\theta_l\right)}{\theta_{l+1}-\theta_l}\left(\theta_n-\theta_{l}\right). 
\end{split}
\label{eq:interpol_1}
\end{equation}
In a second step, the values of \cref{eq:interpol_1} are used and an interpolation is carried out with respect to the accumulated plastic strains. Now, $\alpha_k<\alpha_n<\alpha_{k+1}$ applies and the following holds for the interpolation

\begin{equation}
y\left(\alpha_n,\theta_n\right)=y\left(\alpha_k,\theta_n\right)+\underbrace{\frac{y\left(\alpha_{k+1},\theta_{n}\right)-y\left(\alpha_{k},\theta_n\right)}{\alpha_{k+1}-\alpha_k}}_{=h}\left(\alpha_n-\alpha_k\right).
\label{eq:interpol_2}
\end{equation}
In addition to determining the yield stress $y\left(\alpha_n,\theta_n\right)$, the hardening parameter $h$, which reflects the quotient of yield stress and accumulated plastic strains, is computed simultaneously.
    
\section{Parallel high resolution simulations on modern high-perf\-ormance computers}
\label{sec:hpc}
As described above, discretization in space using finite elements, discretization in time, and using Newton's method to linearize the nonlinear thermo-elastoplastic problem, results in a sparse linear system in each Newton iteration. Caused by the high resolution of our finite element meshes, we have millions of \ac{DOFs} in our computations and thus the solution of extremely large linear systems has to be parallelized using \ac{HPC} resources. In general, iterative solvers combined with efficient preconditioners are well-suited for this purpose. This motivates our choice of the \ac{GMRES} method in combination with suitable domain decomposition preconditioners.

\begin{figure}[]
    \centering
    \includegraphics[clip=true, trim={50 600 50 500},width=1\textwidth]{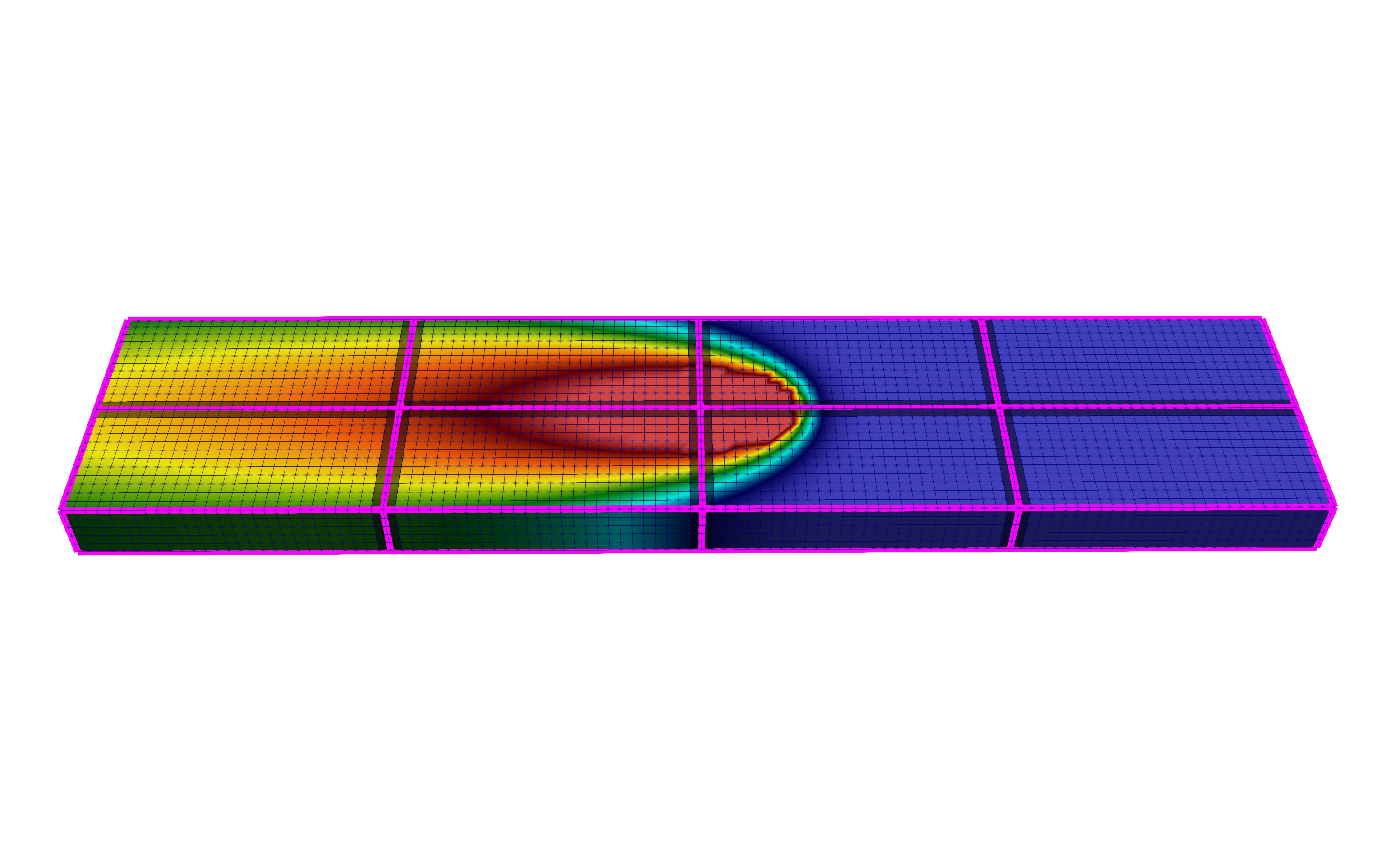}
    \caption{Example of an overlapping domain decomposition of a portion of the thin metal plate into overlapping subdomains. The overlap is visible as the opaque region, while the non-overlapping domain decomposition is highlighted in magenta. Each subdomain is associated to a single MPI rank.} \label{fig:dd_sample}
\end{figure}

\ac{DDMs} are powerful techniques used to solve \ac{PDEs} and the parallelization is based on a divide-and-conquer paradigm. Therefore the computational domain is decomposed into smaller subdomains, the \ac{PDEs} are then localized to the subdomains, and the resulting subdomain problems are solved independently and coupled in an iterative process; see \cref{fig:dd_sample}. Among all overlapping \ac{DDMs}, the  additive Schwarz preconditioner is a popular choice, where information is exchanged within the overlap of neighboring subdomains in each \ac{GMRES} iteration, that is, in each application of the Schwarz preconditioner.  For robustness and scalability, a second level of correction, often referred to as a coarse space, is introduced. This coarse space facilitates the exchange of global information improving the overall convergence by addressing long-range dependencies and ensuring efficient scaling with the number of subdomains. Specifically, in our work, we make use of a two-level Schwarz method with \ac{GDSW} coarse spaces to enhance the efficiency and scalability of the domain decomposition approach; see~\cite{dohrmann2008family,dohrmann2008domain} for details. Let us note that \ac{GDSW} coarse spaces have some important advantages which make them a perfect choice in \ac{LBW} simulations. They can be used with general subdomain shapes, computed algebraically using only information from the system matrix, and are applicable to multi-physics problems, as, for example, thermo-elastoplasticity problems, in a monolithic fashion~\cite{bevilacqua2025highly}. 

The aforementioned Schwarz preconditioners are implemented in the software package FE2TI, a C/C++ library with an MPI/OpenMP~\cite{chandra2001parallel} hybrid parallel implementation exploiting the PETSc library. The FE2TI software has already been shown to scale efficiently to more than one million parallel ranks solving realistic deformation scenarios of dual-phase steels; see~\cite{klawonn2020computational,nakajima} for details. Moreover, as already mentioned, we equipped it with an interface to the FEAP library~\cite{Tay:2008:fea}, to exploit the finite elements and material models from FEAP. Finally, the local and coarse problems in the Schwarz preconditioners are solved exactly using the version of PARDISO implemented in the Math Kernel Library~\cite{intel_onemkl_pardiso,schenk2001pardiso}. 
%To summarize, we use the highly efficient \ac{GMRES} solver provided in PETSc to solve the large and sparse linear systems iteratively, our highly scalable implementation of Schwarz preconditioners with coarse spaces of \ac{GDSW}-type to control the convergence of \ac{GMRES}, and the interface to FEAP to use the material models described in~\cref{sec:FE-model}.

All of these features are used to model the \ac{CTW} test on a portion of the plate. Specifically, we analyze a strip with dimensions $4.44\,\text{mm}$ in width, utilizing highly resolved meshes consisting of millions of \ac{DOFs}, distributed across thousands of parallel compute cores. We simulate $2.4\,\text{s}$ of the \ac{CTW} test, subdivided as follows: from $0\,\text{s}$ to $0.1\,\text{s}$ the heat-up procedure is performed, keeping the laser fixed at its initial position; from $0.1\,\text{s}$ to $1.9\,\text{s}$ the laser moves at a velocity of $1\,\text{m/s}$ (Stage 1 in \cref{sec:ctw_test}); finally from $1.9\,\text{s}$ to $2.4\,\text{s}$, the laser continues to move while an external strain is applied (Stage 2 in \cref{sec:ctw_test}).

Throughout the simulation, boundary conditions on the side faces of the strip are interpolated from the ANSYS simulation of the whole plate, while the laser action is modeled as described in \cref{sec:heat_source}. A time step size of $\Delta t = 0.001\,\text{s}$ is always used (except when explicitly stated), resulting in a total of $\numgru{2400}$ iterations. In each time step, Newton's method is applied until the global residual falls below an absolute tolerance of $10^{-3}$. For the linear systems, the \ac{GMRES} method is employed, which is restarted every 100 iterations, with a stopping criterion of either a relative residual of $10^{-6}$ or an absolute residual of $10^{-10}$ based on the unpreconditioned norm.

\subsection{Convergence study}\label{sec:conv_study}
To choose suitable discretization grid sizes in space and time, we perform a convergence study on the complete \ac{CTW} simulation. We analyze the strip with a width of $4.44\,\text{mm}$, utilizing different, highly-resolved meshes, consisting of $\numgru{324600},\numgru{2330636}, \linebreak \numgru{7573312}$ and $\numgru{17607828}$ \ac{DOFs} distributed across $360$, $720$, $1440$ and again $1440$ cores respectively.

\begin{figure}[t]
    \centering
    \includegraphics[clip=true, trim={150 700 500 800},width=1\textwidth]{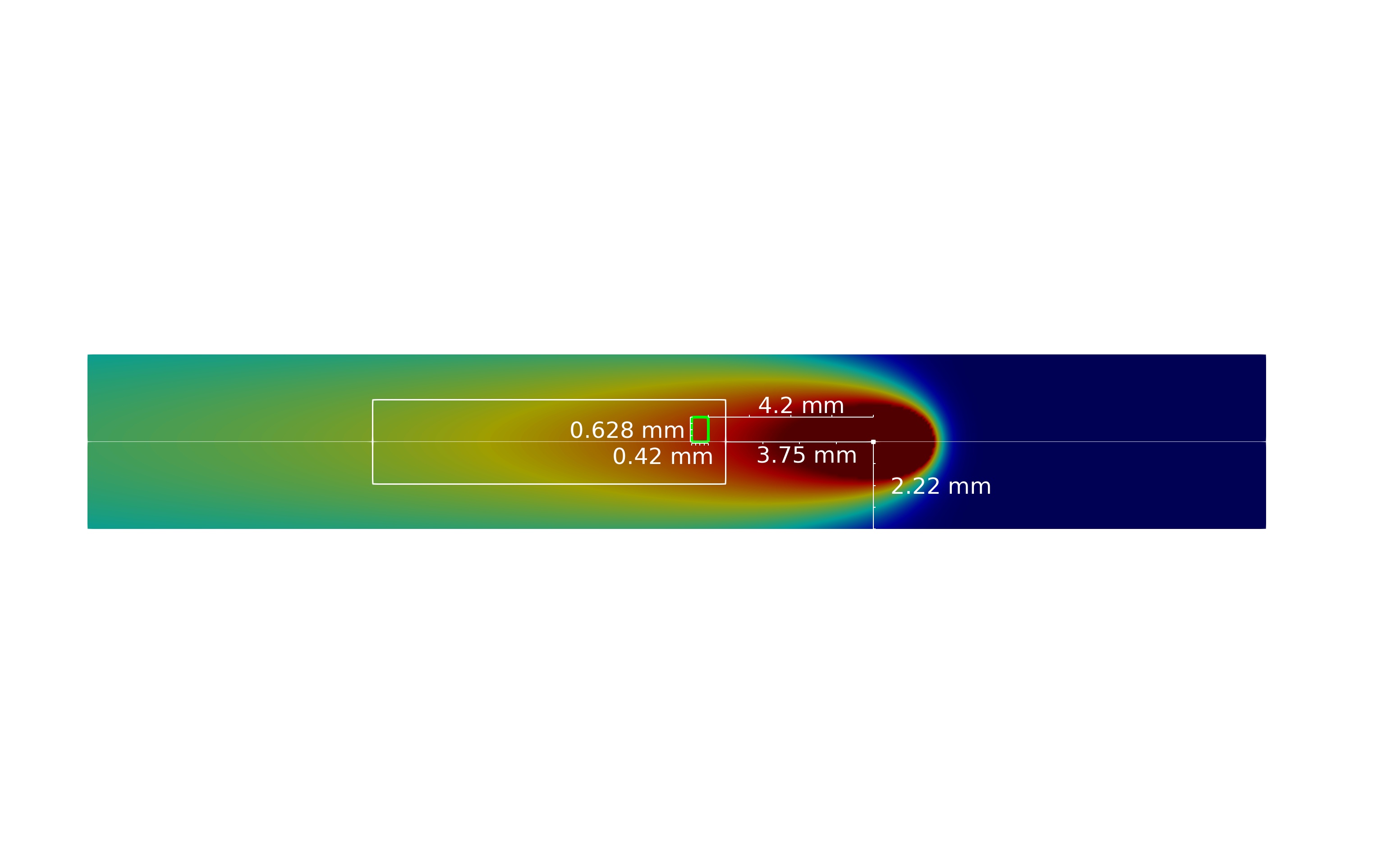}
    \caption{Position of the \ac{ROI} and the small green region over the temperature field.} \label{fig:ROIwind}
\end{figure}

Even if, for the moment, our final goal is to only achieve a qualitative comparisons with the experimental results, a grid convergence study is needed to test our FEM implementation. To this end, we focus our analysis on measuring three quantities of interest located on the surface of the plate regarding the $y$-strain, namely the mean value evaluated within the \ac{ROI} as well as the maximum and mean value evaluated on a small green window (GW) as shown in \cref{fig:ROIwind}. Additionally, since these quantities are computed exclusively on the plate's surface, we also consider the maximum strain value within the \ac{ROI}, extending the analysis along the vertical direction, where we expect to observe higher values. Since we do not perform a quantitative comparison between numerical and experimental results, in all these numerical simulations, the \ac{ROI} size corresponds to the entire \ac{FOV}.

\begin{figure}[]
    \centering
    \includegraphics[width=1\linewidth]{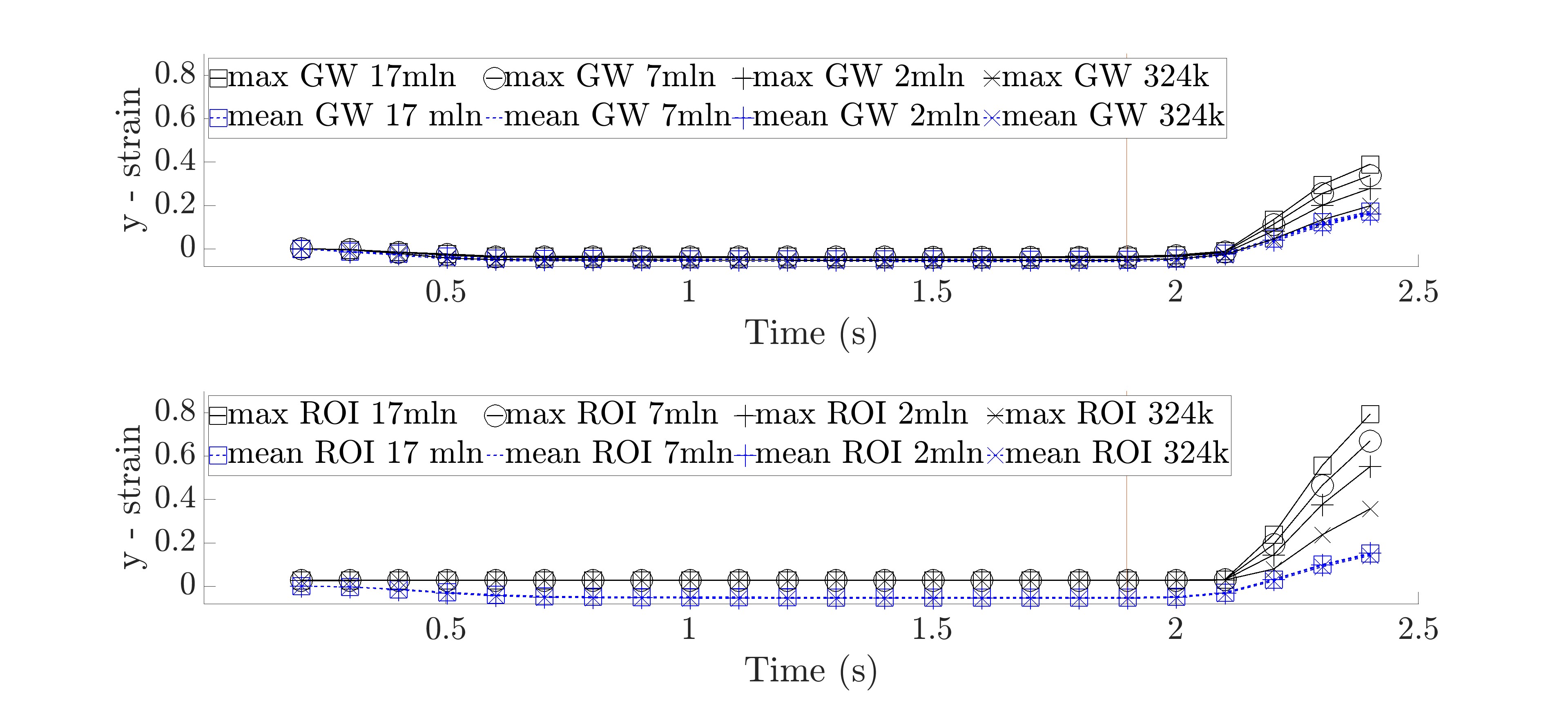}
    \caption{Plot of the quantities of interest over the time steps for the meshes with $\numgru{324600},\numgru{2330636}$, $\numgru{7573312}$ and $\numgru{17607828}$ \ac{DOFs}.} \label{fig:mesh_covergence}
\end{figure}

In \cref{fig:mesh_covergence}, we report these quantities calculate over the time iterations. From these results, we observe that mesh convergence is achieved only up to $2.1\;\mathrm{s}$, which corresponds to the time after $0.2\;\mathrm{s}$ of the application of the external strain. After this time, the plastic effect becomes stronger and the solution seems to be strongly mesh-dependent.

This mesh-dependent material behavior is typical for a stress-strain relation with softening. Based on empirical evidence, we suspect the loss of rank-1 convexity at selected locations; this is also supported by the images of plastic strains in \cref{fig:slicecomp}, which correspond to typical localization bands. Due to the complex thermo-mechanically coupled material behavior and the temperature-dependent material parameters, a more in-depth analysis is required here in the future. To overcome this possibly pathological behavior, existing nonlocal continuum formulations in integral and differential form could be used. The latter, also referred to as gradient-enriched continua, have been successfully applied in the fields of elasticity, plasticity, and damage and provide a robust framework for the analysis of size effects. However, these formulations generally require higher derivatives of certain field functions. Since the qualitative behavior does not change for a finer mesh and in order to save computational resources, for our numerical simulation we use the finer mesh with $\numgru{17607828}$ \ac{DOFs}, corresponding to an element size of about $0.05\,mm$ for the scalability test, and a mesh with $\numgru{7573312}$ \ac{DOFs} for the results in \cref{subsec:CTW}.

% Plot of the yield stranght against the time and equivalent plastic strain on a coarse mesh
%\begin{figure}[]
%    \centering
%    \includegraphics[width=1.1\linewidth]{figures/num_sim/point_eval.jpg}
%    \caption{Top: yield strenght over time, Bottom: yield strenght against the equivalent plastic strain. Mesh with 200000 DOFs. The one with 17 milions is running} %\label{fig:yield}
%\end{figure}

\begin{figure}[]
    \centering
    \includegraphics[width=1\linewidth]{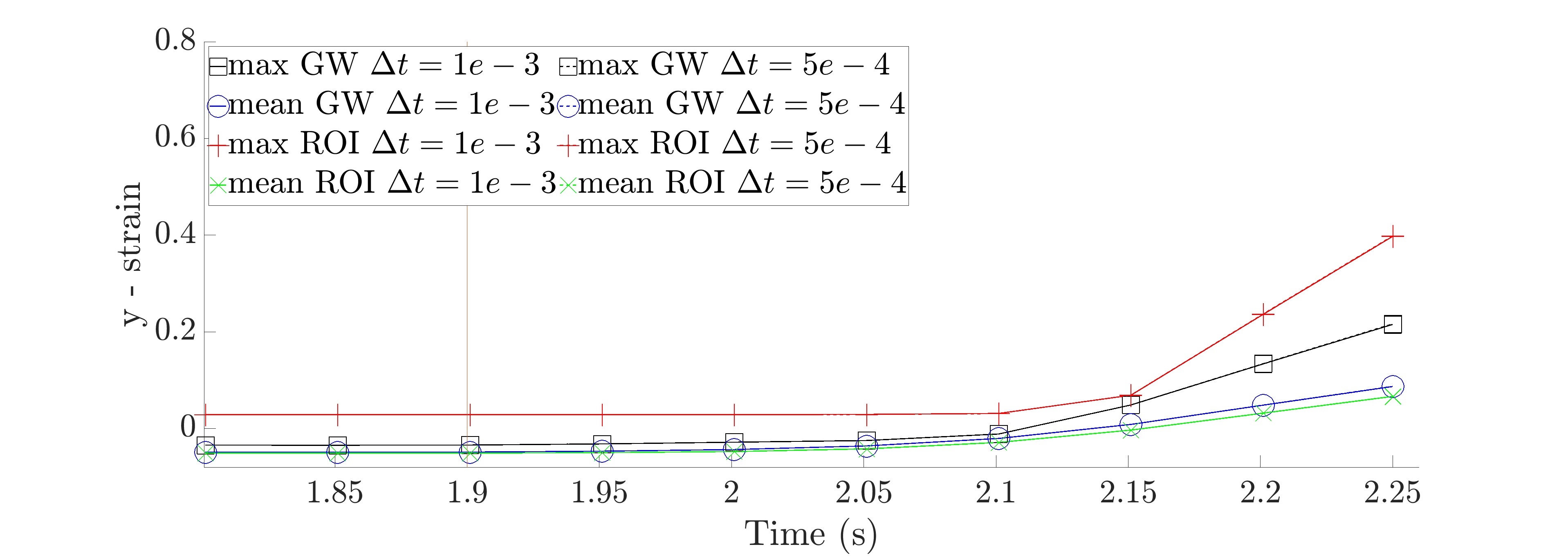}
    \caption{Plot of the quantities of interested for the meshes with $\numgru{17607828}$ \ac{DOFs}, with time steps $0.001\;\mathrm{s}$ and $0.0005\;\mathrm{s}$.} \label{fig:time_covergence}
\end{figure}

To determine the appropriate time step size for the temporal discretization, we conduct a small convergence study. Using the restart option implemented in FE2TI, we simulate the process over the time interval from $1.8\;\mathrm{s}$ to $2.25\;\mathrm{s}$ during which the plastic effects are most pronounced. The study employs the finer mesh with $\numgru{17607828}$ \ac{DOFs}, and the time step size is halved to $\Delta t = 0.0005\;\mathrm{s}$. The results of this analysis, presented in \cref{fig:time_covergence}, show that the convergence in time is reached and a time step of $\Delta t = 0.001\;\mathrm{s}$ is sufficient and therefore used in all following computations.

\subsection{Scalability result}\label{Subsec:scalability}

\begin{figure}[t]
    \centering
    \begin{tabular}{cc}
        \raisebox{1.1\height}{
        \begin{tabular}{rrrrrr}
        \hline
        $Cores$ & $Coarse$  & $T_{PC}$ & $T_{GMRES}$ & $T_{Ass}$ & $T_{Tot}$ \\
        \hline\hline
        90   & 1062  & 234.8 & 37.0 & 33.0 & 313.6 \\
        180  & 2142  & 104.2 & 24.4 & 17.0 & 150.6 \\
        360  & 5358  & 51.6  & 12.2 & 9.2  & 82.0  \\
        720  & 11783 & 18.4  & 6.9  & 5.2  & 33.4  \\
        1440 & 23656 & 10.0  & 4.8  & 3.8  & 21.0  \\
        2880 & 49468 & 6.04  & 6.36 & 1.9  & 15.2  \\
        \hline
        \end{tabular}   
        }
        &
        \hspace{-0.3cm}\includegraphics[clip=true, trim={10 20 110 20},width=0.3\textwidth]{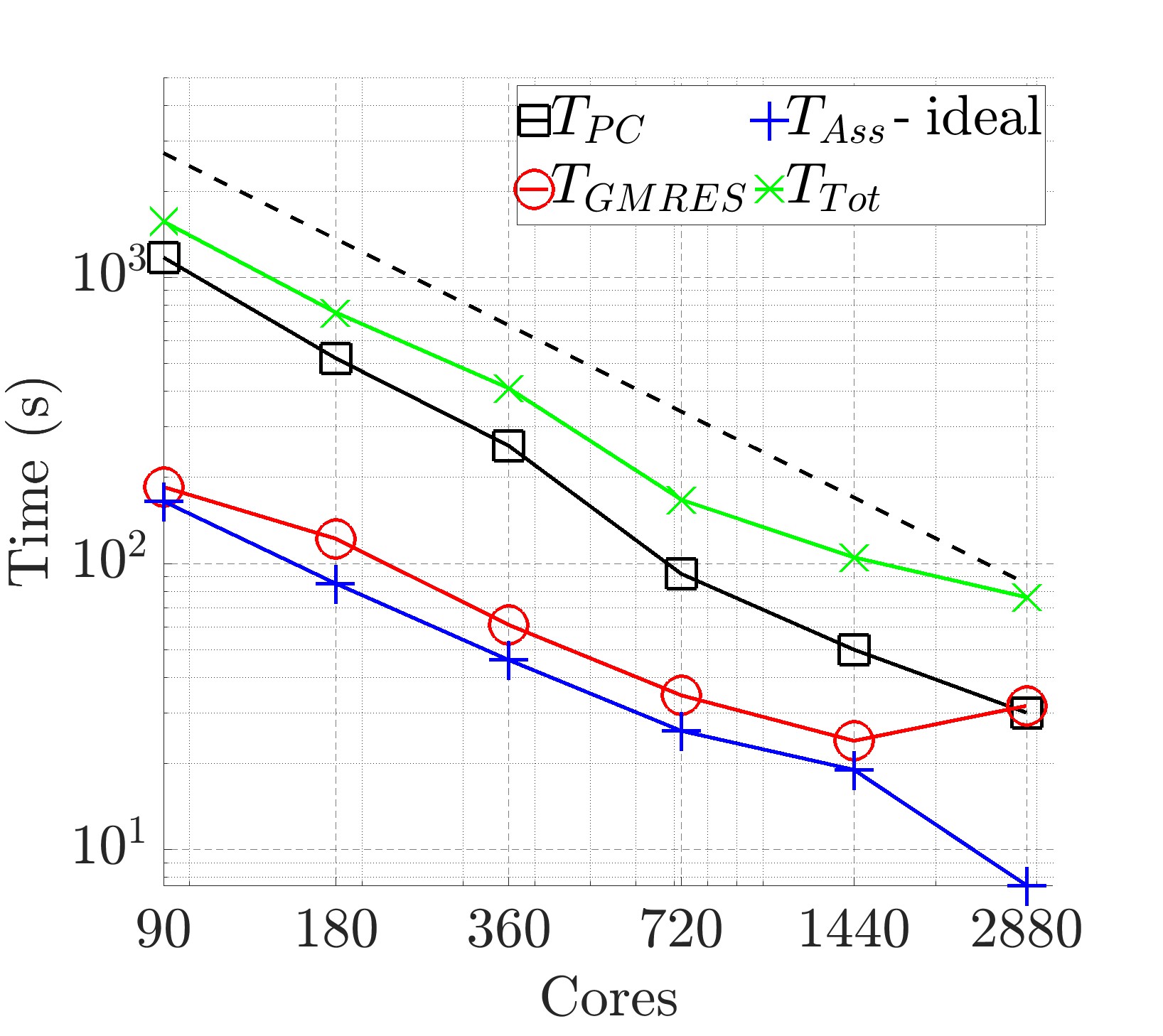}        
    \end{tabular}
    \caption{ Strong scalability test for the first $0.005\,\mathrm{s}$ of a LBW process simulation with $\Delta t = 0.001\,\mathrm{s}$ and $\numgru{17607828}$ million \ac{DOFs}. On the left, the table reports the average timings per time iteration, where $Cores$ represents the number of MPI ranks, $Coarse$ is the size of the coarse space, $T_{PC}$ is the time to construct the preconditioner, $T_{GMRES}$ is the time to compute the solution, $T_{Ass}$ is the time to assemble the FEM matrices and $T_{Tot}$ is the total time per iteration. On the right, the logarithmic plot shows the cumulative timings over the LBW simulation.}\label{fig:StrongScaling}
\end{figure}

To prove the efficiency of our parallel software for the explicit application of \ac{LBW} we perform a parallel strong scaling test. Therefore, we choose a setup with a fixed mesh resolution and measure the acceleration regarding the time to solution when increasing the number of computational cores.  More precisely, we use a setup close to a production run, that is, we fix the total problem size to $\numgru{17607828}$ million \ac{DOFs} while scaling up the resources. In \cref{fig:StrongScaling}, we report the average computing times required to assemble the problem, build the preconditioner, and solve the linear systems.

 In this strong scaling test, and in all the other numerical simulations in this work, we are using a two-level additive Schwarz preconditioner with an overlap of one layer of finite elements and a \ac{GDSW}-R\ac{GDSW} coarse space designed for a thermo-elastoplasticity problem; see~\cite{bevilacqua2025highly} for more details. To enhance the efficiency of our solver, we employ a recycling strategy and the preconditioner is only set up once in each time step; see~\cite{bevilacqua2025highly} for details. 

The results generally indicate a good scaling behavior and the runtime decreases as additional cores are utilized. While scalability is nearly perfect up to 720 cores, we observed  less improvement when scaling from 720 to \numgru{1440} cores and a further deterioration from \numgru{1440} to \numgru{2880}. Regarding the preconditioner timings, this is primarily due to the direct solver required to handle the coarse space, which becomes large at this scale while, at the same time, the local subdomain problems become smaller and smaller. Let us note that for a larger problem with more than 17 million \ac{DOFs} also more than 1\,440 cores can be used efficiently. All in all, we are able to reduce the average time for one time step from $313.3\;\mathrm{s}$ on 90 cores to $21.0\;\mathrm{s}$ exploiting \numgru{1440} cores, which corresponds to an excellent parallel efficiency of 93\%. 
Due to the motivations explained before, we note that the parallel efficiency drops to $64\%$ when computed on \numgru{2880} cores. Therefore, for this specific problem setup, the range between 720 and \numgru{1440} cores seems to be optimal.

\subsection{Numerical simulation of \ac{CTW} test and comparison with real experiment}\label{subsec:CTW}

 In order to compare our numerical simulations with the physical experiments, we utilize a high-resolution mesh consisting of $\numgru{7573312}$ \ac{DOFs}, distributed across 1440 cores. Again, we simulate $2.4\,\text{s}$ of the \ac{CTW} test, with a time step size of $\Delta t = 0.001\,\text{s}$. We report that the complete simulation requires, on average, approximatively 38 GMRES iterations for each linear system with 8 Newton iterations in each time step. The total time for the simulation is $49\,875\;\mathrm{s}$ that corresponds to an excellent $20.8\;\mathrm{s}$ per time step. Let us note that we chose the parameters, especially the time step size, the mesh resolution, the number of cores, and the preconditioner setup, based on the convergence tests and scalability tests documented in \cref{sec:conv_study} and \cref{Subsec:scalability}.

\begin{figure}[]
    \centering
        \includegraphics[clip=true, trim={100 200 50 200},width=1\textwidth]{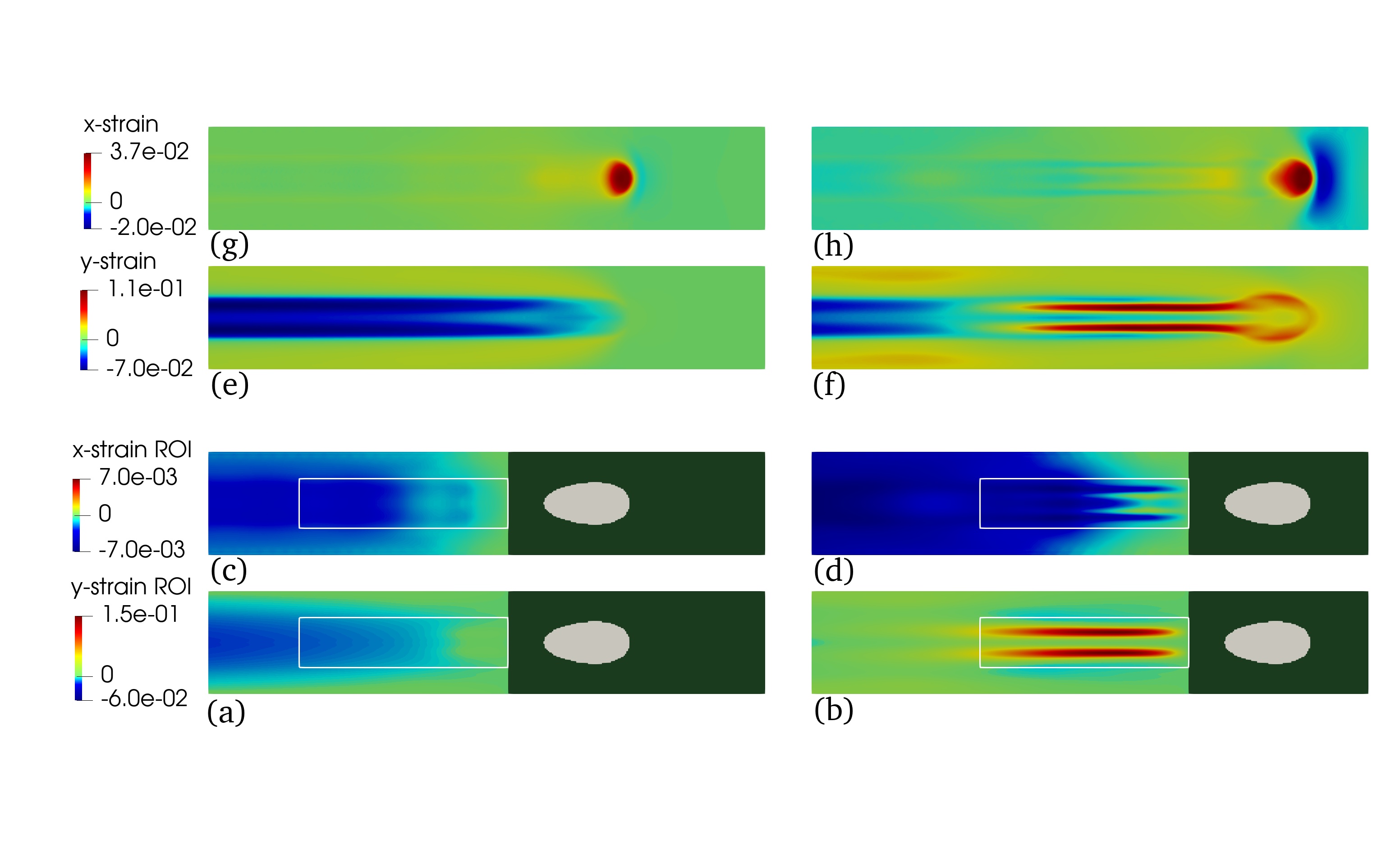} \\
    \caption{Strain distribution on the surface of a fixed section of the plate from $24$ to $44$ mm. From top to bottom: the $x$-strain, $y$-strain, and again the $x$- and $y$-strain with the optical framework technique. The white rectangle shows the position of the \ac{ROI} with respect to the melting pool. On the left columns the distribution of these values is represented at $2.0\;\mathrm{s}$, while on the right one at $2.2\;\mathrm{s}$.}
    \label{fig:ystrain}
\end{figure}

In \cref{fig:ystrain}, we show the distribution of different strain evaluations, at the time points $t=2.0\;\mathrm{s}$ (left column) and $t=2{.}2\;\mathrm{s}$ (right column), computed during the complete simulation. From top to bottom, the first two rows report the total strains in the $x$- and $y$-direction, that is the strains in the welding and its orthogonal direction, respectively, while the third and fourth rows depict the same quantities computed by using an approach similar to the optimal framework technique. This means that in the two last rows, we only depict the fluctuations after a point is covered by a moving region. This area, that includes the \ac{ROI}, is placed $3{.}75\;\mathrm{mm}$ behind the welding seam, and a reintialization of the displacements to the current value takes place.

The results show that higher values for the total strains in the $x$-direction occur in the area of the melt pool. This correlates with the high temperatures and the tendency of the material to expand. 
Outside the melt pool, little to no strains appear while slight tensile strains develop in front of the melt pool (\cref{fig:ystrain}g,h). For the strains in the $y$-direction, the data indicates that compressive strains initially form in and behind the weld pool, while the surrounding area exhibits little to no strain (\cref{fig:ystrain}e). Over time and with increasing external load, the distribution changes, and the applied tensile strains affect the region directly behind the melt pool. Two strips of tensile strains become visible (\cref{fig:ystrain}f). 
For the fluctuations of the strains in $x$-direction, the deviations range between one fifth and one third of the total strains. 
However, increased pressure zones develop in the area behind the melt pool due to the cooling that has already begun and the further heating of the area in front (\cref{fig:ystrain}c,d). The fluctuations in $y$-direction reach the same order of magnitude as the total strains (\cref{fig:ystrain}a,b). The external load likely plays a primary role in this behavior. The two strips of tensile strains at time $t=2{.}2\;\mathrm{s}$ also appear (\cref{fig:ystrain}b). The absolute values compared to the total strains increase for the tensile range. This results from the fact that before entering the \ac{ROI}, a compression range existed, which has now shifted to a tensile range. As a result, the  difference between the values can exceed the absolute value.

%\begin{figure}[]
%    \centering
 %       \begin{tabular}{c}
 %                         \includegraphics[width=0.7\textwidth]{figures/1700_x.PNG} \\
  %          \hspace{-2mm} \includegraphics[width=0.7\textwidth]{figures/1700.PNG} \\
   %     \end{tabular}
    %\caption{This figure shows the distribution of the x and y strain short before external load application in observed ROI \AG{ I must check a sign of x- strain. This seems to be wrong one.} {\PH Not yet mentioned, as there is no direct comparison with the simulations. In addition, the x-values should still be incorrect.} }
  %  \label{fig:exp_1700}
%\end{figure}

\begin{figure}[]
    \centering
    % \begin{subfigure}[t]{\textwidth}
        \centering
        \includegraphics[width=0.7\textwidth]{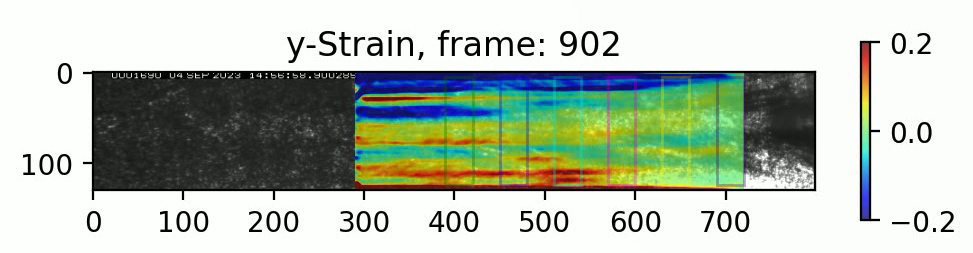}
        \includegraphics[width=0.7\textwidth]{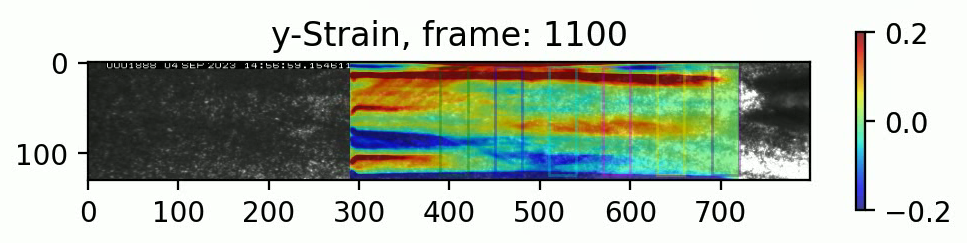}
        % \caption{y-strains for  frames 902 and 1100}
        % \label{fig:exp_1904_a}
    % \end{subfigure} \\
    % \begin{subfigure}[t]{\textwidth}
    %     \centering
    %     \includegraphics[width=0.7\textwidth]{figures/CTW-20-6_1-1_2_Ver8 Strain graph.png}
    %     \includegraphics[width=0.7\textwidth]{figures/CTW-20-6_1-1_2_Ver8 Displacement graph.png}        \caption{{\PH \sout{$x$- and }}$y$-strain {\PH and $y$-displacement} components over the frames}
    %     \label{fig:exp_1904_b}
    % \end{subfigure}
    \caption{Distribution of the $y$-strain shortly before (frame 902)  and shortly after (frame 1100) the external load application in the observed \ac{ROI}}
    \label{fig:exp_1904}
\end{figure}

\cref{fig:exp_1904} shows the experimental results in comparison to the numerical results. The diagram illustrates the $y$-strain shortly before (frame 902) and shortly after (frame 1100) the start of the application of the external load. The frame rate of the recordings corresponds to 778 fps, see \cref{sec:opt_meas_tech}. Thus, one frame equals approximately $1{.}285*10^{-3}$ seconds, while frame 1000 corresponds approximately to the time of load application. Although the absolute values of the strains are not identical to the numerical results, the strip-shaped profile of a tension band behind the melt pool on the outer area of the previously melted area can also be seen here directly after the application of the external load. The reasons for the differences in the absolute values of the strains can be manifold. On the one hand, the optical framework technique employs a relatively coarse meshing, which, in comparison to the highly refined FEM mesh, may lead to less precise results. Furthermore, in an experimental setting, no simplifications can be made, whereas the numerical simulation only considers a thermo-elastoplastic material behavior on the macroscale, while the model completely disregards microstructural processes occurring during cooling. 

\begin{figure}[]
    \centering
        \begin{tabular}{c}
            \includegraphics[clip=true, trim={100 300 0 650},width=0.9\textwidth]{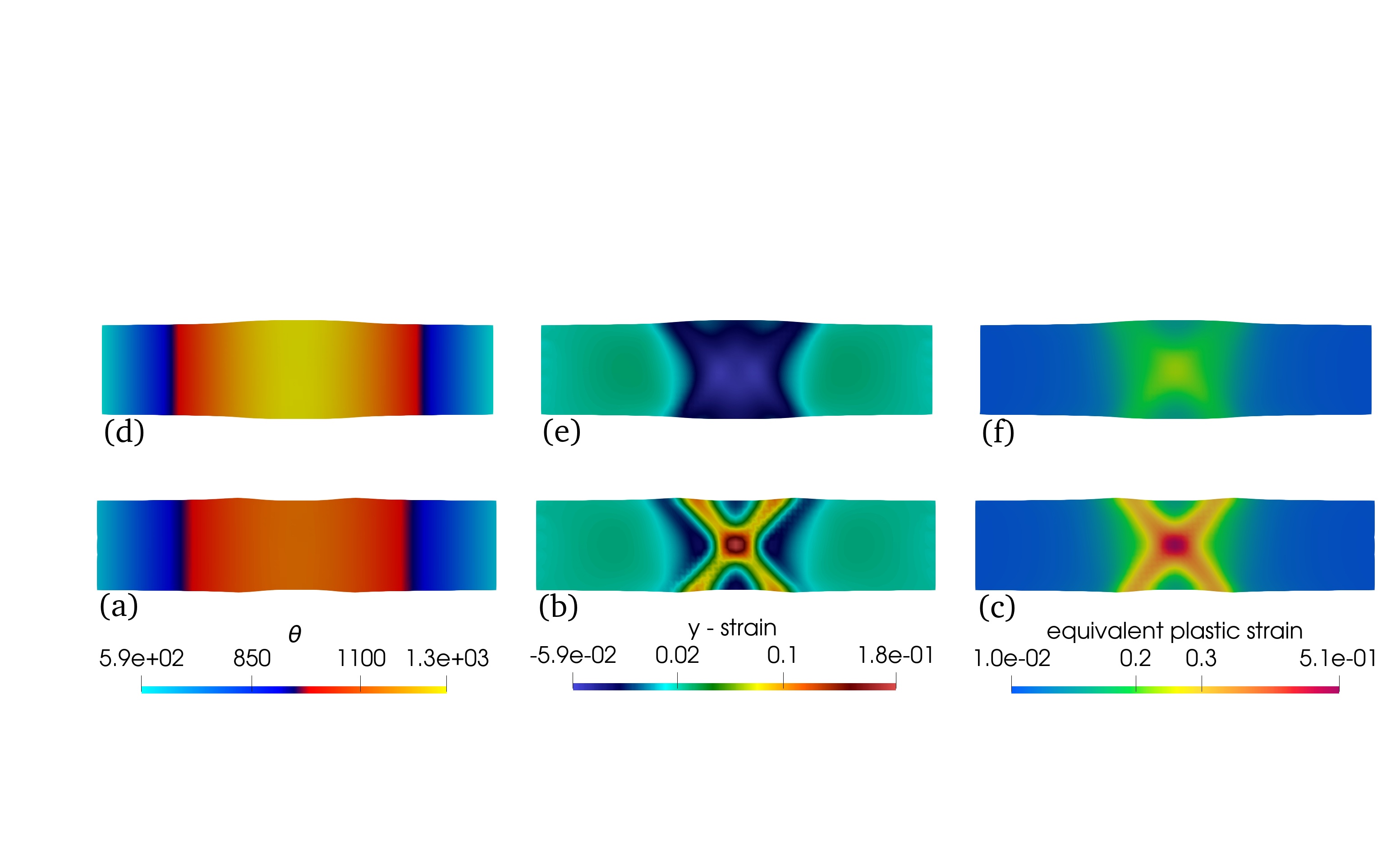} \\
        \end{tabular}
    \caption{From left to right: temperature $\theta$, $y$-strain and equivalent plastic strain of a cross section of the plate positioned $35\,\mathrm{mm}$ from the edge. The slice are taken at two different time steps, namely $2.1\,\mathrm{s}$ and $2.2\,\mathrm{s}$.}
    \label{fig:slicecomp}
\end{figure}

\cref{fig:slicecomp} depicts the cross-section at a distance of $35\;\mathrm{mm}$ from the edge where the welding process began, shown for the time points $t=2{.}1\;\mathrm{s}$ and $t=2{.}2\;\mathrm{s}$. The figure presents the temperature, the strains in $y$-direction, and the equivalent plastic strains. The temperature field exhibits a uniformly distributed temperature within the solidifying material for both time points. However, in the strain distribution, a cross-shaped pattern emerges, becoming more visible as time progresses. This cross-shaped distribution indicates localization bands. This means that the strains increase significantly and become localized. The underlying cause is typically the interaction between thermal and mechanical effects; see \cite{Mie:1995:ato}.

Based on the localization, we can conclude that the cracks observed in the experiment originate from inside the workpiece and later propagate to the surface. This aspect can not be identified by the optical frame technique, since this approach only detects changes at the surface and is therefore unable to capture cracks that emerge from within the material.

%--- Section ---%
\section{Conclusions}\label{sec:conclusion}
In the present work, large-scale thermo-mechanical simulations of LBW processes, more precisely simulations of the \ac{CTW} test, were conducted using highly resolving discretizations and high-performance computing to obtain high-fidelity numerical results. The chosen hybrid approach combining the strengths of the commercial software ANSYS and the high-performance computing framework FE2TI resulted in detailed virtual insights into the welding process. Especially the divide-and-conquer parallelization paradigm used in FE2TI enabled large-scale simulations exploiting modern supercomputers and led to a very fine resolution of the surrounding of the weld seam. 
The numerical results are in good qualitative agreement with the experimental ones, particularly with respect to strain distribution and localization phenomena, despite the complexity of the thermo-mechanical coupling in LBW. However, the quantitative agreement is still limited, which may be due to various factors, such as not considering the dendritic microstructure in the simulation, the coarse meshing of the optical framework in the experiment, or the non-consideration of various processes in the melt pool. Therefore, to improve the predictive capabilities of numerical simulations and to further investigate the influence of microstructural effects on crack initiation, future work should focus on multi-scale modeling approaches for the simulations.
Furthermore, to fully trust the FE simulation, the results should converge using increasingly finer resolution and more DOFs. This convergence has been observed for the mean strain, it also holds for the maximum strain until the plastic effect becomes dominant. The latter indicates that further localization analysis is required in the future to better capture strain concentration effects. It also highlights the mesh dependency of the solution, with finer meshes leading to a softer material response.
To summarize, exploiting the potential of \ac{HPC} resulted in highly resolved, coupled thermo-mechanical simulations of LBW processes and a good qualitative agreement with experiments. For a quantitative prediction of solidification cracks, further steps must be taken, such as the consideration of the multiscale nature of this process. On the other hand, the experimental evaluation algorithms must be further improved, to minimize the accumulated error and provide a higher accuracy especially for strain computations based on the obtained displacement distribution. This must involve {a further improvement of the image resolution as well as the application of different numerical filters.      

%-------------------------------------------
% Optional Contents
%-------------------------------------------

%--- Section ---%
%PH: Apparently not necessary for arxiv:
% \section*{Conflicts of Interest} 
% Declare conflicts of interest or state “The authors declare no conflict of interest.” Authors must identify and declare any personal circumstances or interests that may be perceived as inappropriately influencing the representation or interpretation of reported research results. A detailed definition of conflicts of interest is available at the following site: \url{https://academic.oup.com/journals/pages/authors/preparing_your_manuscript/ethics#conflict}.

%--- Section ---%
%PH: Apparently not necessary for arxiv:
% \section*{Author Contributions}
% Must include all authors, identified by initials, for example: “Conceptualization, S.R.. and D.A.; methodology, S.R..; software, S.R..; validation, S.R.., D.A. and K.L.; formal analysis, S.R..; investigation, S.R..; resources, S.R..; data curation, S.R..; writing—original draft preparation, S.R..; writing—review and editing, S.R..; visualization, S.R..; supervision, S.R..; project administration, S.R..; funding acquisition, D.A.” Individual contributions are specified according to NISO CrediT (Contributor Roles Taxonomy) described at the following site: \url{https://credit.niso.org/}.

%--- Section ---%
\section*{Acknowledgments}
This research is funded by the Deutsche Forschungsgemeinschaft (DFG, German Research Foundation) – 434946896 within the research unit FOR 5134 ``Solidification Cracks during Laser Beam Welding: High Performance Computing for High Performance Processing''. The authors gratefully acknowledge the scientific support and HPC resources provided by the Erlangen National High Performance Computing Center (NHR@FAU) of the Friedrich-Alexander-Universit\"at Erlangen-N\"urnberg (FAU) under the NHR project k109be10. NHR funding is provided by federal and Bavarian state authorities. NHR@FAU hardware is partially funded by the German Research Foundation (DFG) - 440719683. 
Interchange and data storage was carried out by exploiting Kadi4Mat~\cite{brandt2021kadi4mat}.

%-------------------------------------------
% References
%-------------------------------------------

%-------------------------------------------
% Appendix
%----------------------------------------

\bibliographystyle{siamplain}
\bibliography{literature}
\end{document}